\documentclass{article}
\usepackage{fullpage}
\usepackage{url,mathtools,bm}

\usepackage{enumerate,amsfonts, graphicx, amsmath}
\usepackage{graphicx}
\usepackage{epsfig}
\usepackage{multirow,color}
\newcommand\llbracket{[\kern-2pt[}
\newcommand\rrbracket{]\kern-2pt]}
\newcommand\argmin{\mathop{\rm argmin}}

\renewcommand\a{{\bm{a}}}
\renewcommand\b{{\bm{b}}}
\renewcommand\c{{\bm{c}}}
\renewcommand\d{{\bm{d}}}
\newcommand\e{{\bm{e}}}
\newcommand\f{{\bm{f}}}
\newcommand\g{{\bm{g}}}

\newcommand\R{{\bm{R}}}
\newcommand\s{{\bm{s}}}

\renewcommand\u{{\bm{u}}}

\renewcommand\v{{\bm{v}}}

\newcommand\x{{\bm{x}}}

\newcommand\blambda{{\bm{\lambda}}}
\newcommand\bxi{{\bm{\xi}}}

\newcommand\bz{{\bm{0}}}

\newcommand{\trace}{\mathop {\rm Tr}}

\title{Second-order cone
  interior-point method for quasistatic and moderate dynamic cohesive fracture\footnote{Supported
in part by an NSERC Discovery Grant.}}
\author{Stephen A. Vavasis\thanks{Department of Combinatorics \&
Optimization, University of Waterloo, 
200 University Ave.~W., Waterloo, Ontario, N2L 3G1, 
Canada, vavasis@uwaterloo.ca.} \and
Katerina D. Papoulia\thanks{Department of Applied Mathematics
and Department of Mechanical \& Mechatronics
Engineering,
University
of Waterloo, 
200 University Ave.~W., Waterloo, Ontario, N2L 3G1,  
Canada, papoulia@uwaterloo.ca} \and
M.~Reza Hirmand\thanks{Department of Mechanical \& Mechatronics
Engineering, 
University
of Waterloo, 
200 University Ave.~W., Waterloo, Ontario, N2L 3G1,  
Canada, mhirmand@uwaterloo.ca}}
\begin{document}
\maketitle
\begin{abstract}
  Cohesive fracture is among the few techniques able to
  model complex fracture nucleation and propagation
  with a sharp (nonsmeared) representation
  of the crack.  Implicit time-stepping schemes are often favored 
  in mechanics due to their ability to take larger time steps in
  quasistatic and moderate dynamic problems.  Furthermore,
  initially rigid cohesive models are typically preferred when
  the location of the crack is not known in advance, since
  initially elastic models artificially lower the material stiffness.
  It is challenging to include an initially rigid
  cohesive model in an implicit scheme because
  the initiation of fracture corresponds
  to a nondifferentiability of the underlying potential.  In
  this work, an interior-point method is proposed for implicit time
  stepping of initially rigid cohesive
  fracture.  It uses techniques developed for convex second-order
  cone programming for the nonconvex problem at hand.  The underlying cohesive model
  is taken from Papoulia (2017) and is based on a nondifferentiable
  energy function.  That previous work proposed an algorithm based on successive
  smooth approximations to the nondifferential objective for solving
  the resulting optimization problem.  It is argued herein that cone
  programming can capture the nondifferentiability without smoothing,
  and the resulting cone formulation is amenable to interior-point
  algorithms.  A further benefit of the formulation is that other
  conic inequality constraints are straightforward to incorporate.
  Computational results are provided showing that certain contact
  constraints can be easily handled and that the
  method is practical.
\end{abstract}
\maketitle

\section{A nondifferentiable energy model for cohesive fracture}
In this section, we review prior work by Papoulia \cite{optiIJF}
as well as related
works
on nondifferentiable energy models for fracture.  We assume isothermal
conditions throughout so that
temperature effects can be neglected
and consideration of thermodynamics reduces to mechanical
potential and kinetic energy.

Let $\Omega\subset\R^{n_{\rm dim}}$ ($n_{\rm dim}=2,3$) denote the initial configuration
of the body under consideration.  Let $\mathcal{S}\subset\Omega$ be a union
of $(n_{\rm dim}-1)$-dimensional surfaces (when $n_{\rm dim}=3$) or
curves (when $n_{\rm dim}=2$) that may each cut across the entire domain.
Let
$u:\Omega\rightarrow\R^{n_{\rm dim}}$ be the displacement field, assumed
to be differentiable except possibly for jumps on $\mathcal{S}$.

In this model of solid mechanics and fracture, two
potential energies exist, one associated with the bulk model and one
with a network of surfaces $\mathcal{S}$ inside the domain that serve as
potential sites of fracture.  Thus, the mechanical potential
has the form:
\begin{equation}
 \mathcal{E} = 
 \int_\Omega \Psi(u)\,dV + \int_{\mathcal{S}} \Phi(\llbracket u\rrbracket)\,dS
 - \int_\Omega f\cdot u\,dV
 - \int_{\mathcal{S}_t}\tilde{t}\cdot u\,dS
  \label{eq:twoterm} 
 \end{equation}
where $\llbracket u\rrbracket$ denotes the jump in $u$ across the
surface.  Here, $\Psi$ corresponds to the strain-energy density function
while $\Phi$ corresponds to energy density of new surfaces (fracture).
The final two terms correspond to body loads ($f$) and traction loads ($\tilde t$)
respectively, Here $\mathcal{S}_t$ denotes the portion of $\partial\Omega$
with traction loads.
Note that the use of a potential energy functional implies reversibility;
we return to this matter below.
Later on, we will add another term to the optimization formulation to
account for momentum. In addition, we will impose displacement
boundary conditions
and inequality constraints.  The latter will be used to model contact
and a no-interpenetration requirement for the mesh.

We further stipulate that $\Phi(\llbracket u\rrbracket)$
is a nondifferentiable function of $\llbracket u\rrbracket$ 
when $\llbracket u\rrbracket =0$.  As explained in
\cite{optiIJF}, this is an essential ingredient of the formulation; see
also Charlotte et al.\ \cite{charlotte} in a slightly different context.
Because of the nondifferentiability at $\llbracket u\rrbracket =0$, no jumps in $u$
will occur across any surface
until a positive finite level of loading occurs.

The above model falls into the category of ``cohesive zone'' models
\cite{Barenblatt,Rice,Dugdale} because it accounts for crack propagation
with explicit representation of crack surfaces and an associated
displacement-traction relation (which is obtained as a derivative of $\Phi$ for
nonzero values of $\llbracket u\rrbracket$).
Furthermore, it falls into the category of ``initially rigid'' cohesive
zone models because of the property that there is no crack opening until
a specific positive finite load level is attained.   Initially rigid models are
preferred over the alternative ``initially elastic'' models in problems
where the crack path is not known a priori.  Inclusion of a network
of initially elastic
surfaces would lower the global stiffness of $\Omega$; as the number of
crosscutting surfaces in $\mathcal{S}$ tends to infinity, the global stiffness is driven to 0.
In contrast, there is no limit to how much surface area may be encompassed
by $\mathcal{S}$ in \eqref{eq:twoterm} for
the class of nondifferentiable potentials $\Phi(\llbracket u\rrbracket)$
proposed in \cite{optiIJF}.

 In \cite{PV1}, it was argued that unless
significant care is taken in designing the algorithm, methods for
initially rigid cohesive fracture are
likely to be ``time discontinuous.''
The issue is
that after space discretization, a system of ODE's for nodal values of the
displacement $u$ and other quantities arises, i.e., a system of the form
$d\u/dt = \f(\u)$.  The forcing
function $\f(\u)$ of these ODE's is a discontinuous
function of $\u$, and this leads to nonconvergent or unreliable numerical methods.
In \cite{optiIJF} and also in this paper,
the problems of time discontinuity are sidestepped
because the modeling technique does not lead to a system of ODE's---the
usual step of passing to a weak form does not apply because
the potential is nondifferentiable. Instead, the method involves time
steps each of which corresponds to a
physically based energy minimization operation.

The formulation \eqref{eq:twoterm} thus reduces the problem of modeling fracture
to a sequence of optimization problems.  These are infinite dimensional 
problems,
but they are reduced to finite-dimensional optimization using 
finite element analysis as discussed in Section~\ref{sec:fem}.
This problem was solved in \cite{optiIJF}
using a continuation method.  Hirmand and Papoulia \cite{DGopti} solve
it using a Nitsche discontinuous Galerkin method
(see also the related work by Radovitzky et al.\ \cite{Radov})
in which the multipliers of the
optimization problem are interpreted as stresses at the crack surface.
The contribution of the present paper is a solution method for the optimization problem
\eqref{eq:twoterm} using a novel interior-point method.  A
key step in the development, as
explained in Section~\ref{sec:ipfm},
is to recast a certain equation (namely, the second
line of \eqref{eq:optmodeldisc}) that appears in the optimization problem
as an inequality (namely, the second line of \eqref{eq:mainprob}).

This technique is commonplace in the optimization literature
but is new (as far as we know) to fracture mechanics.
General background on interior-point methods is provided in Section~\ref{sec:ipbg}.
As mentioned earlier, the development of the formulation
continues in Section~\ref{sec:fem}, which explains our
finite-element discretization.
The method as described so far is reversible.
Irreversibility may be incorporated via the
additional dependence of $\Phi$ in \eqref{eq:twoterm} on
a damage variable as detailed in Section~\ref{sec:fem}.

The interior-point
formulation is provided in Section~\ref{sec:ipfm}.  Most
of the literature on interior-point methods relates to convex
optimization.  Our optimization problem is nonconvex, which
requires modifications to the interior-point method compared
to the previous literature as explained in Section~\ref{sec:nonconv}.
The interior-point method needs a feasible starting point; for
this we rely on a technique developed in Section~\ref{sec:phaseI}.
Details of the computational procedure are spelled
out in Section~\ref{sec:compproc}.  Our computational
experiments are described in Section~\ref{sec:compexp}; these
experiments involve checking the balance of energy, the computation of
which is described in Section~\ref{sec:energy_balance}.  We conclude
in Section~\ref{sec:conc} with an itemization of the development
of the optimization models as well as the components of our computational
method.

We conclude this section with a discussion of related literature.
Other
than Papoulia \cite{optiIJF} and Hirmand and Papoulia \cite{DGopti},
the most closely related work is Lorentz's \cite{Lorentz2008subgradient} method, which
also treats initially rigid fracture using a potential like
\eqref{eq:twoterm} for the same reasons as us.  Lorentz does not use
an optimization method per se but rather considers the subdifferential
of \eqref{eq:twoterm} as a generalization of a system of equations for
generating a time step.
Slightly more distantly related to the present work
is the phase-field method of modeling
fracture \cite{bourdin2008}.
In this case, energy minimization is also invoked, but the
functional pertains to a smeared crack location rather than a sharp
surface.  As a consequence, a sharp representation of the crack 
must be determined {\em a posteriori}, although some authors
e.g., Geelen et al.\ \cite{Dolbow2},
Wang and Waisman \cite{Waisman},
have shown recently
that a sharp representation of the crack can be
directly coupled to a phase-field model.

\section{Interior-point algorithms}
\label{sec:ipbg}

In this section we present general background on interior-point
methods.  For more in-depth treatment, see, e.g., \cite{Ye}.
The application of these methods to initially rigid
cohesive fracture is provided in Section~\ref{sec:ipfm}.

A {\em closed convex cone} is defined to be a set $K\subset \R^n$
with the properties that (i) $K$ is closed, (ii) $K$ is convex,
(iii) $K$ is a cone, i.e., $\x\in K\Rightarrow\lambda\x\in K$ for
all $\lambda \ge 0$, and (iv) $\bz\in K$ (i.e., $K\ne\emptyset$).

Two important special examples of closed convex cones are 
$\R_n^+=\{\x\in\R^n: x_i\ge 0\>\forall i=1,\ldots,n\}$, the
{\em nonnegative orthant},  and
$C_2^{n}=\{\x\in\R^n: x_1\ge \Vert \x(2:n)\Vert\},$ the {\em second-order
  cone}.  Here, $\x(2:n)$ (Matlab notation) denotes the subvector of $\x$
indexed by coordinates $2$ through $n$.

The two cones mentioned in the previous paragraph both have standard
{\em self-concordant barrier functions}.  For 
$\R_n^+$, the standard self-concordant barrier function is 
$\phi_{\rm NNO}(\x)=-\sum_{i=1}^n\log(x_i)$.  For $C_2^{n}$, the
standard self-concordant barrier function is
$\phi_{\rm SOC}(\x)=-\frac{1}{2}\log(x_1^2-x_2^2-\cdots-x_n^2)$.  We regard
these functions as taking on the value infinity outside the relevant cones.
These functions have the property that they are strictly convex functions
on the interior of their respective cones, and they tend to
infinity as the boundary of the cone is approached.  ``Self concordance''
involves two other technical properties; see \cite{NemNest}.

If $K_1\subset \R^{n_1}, \ldots,K_r\subset \R^{n_r}$ 
are all closed convex cones with barrier functions
$\phi_1,\ldots,\phi_r$,
then $K_1\times\cdots\times K_r$ is also a closed convex cone, and its 
barrier function is $\phi_1(\x_1)+\cdots+\phi_r(\x_r)$ for 
$(\x_1,\ldots,\x_r)\in K_1\times\cdots\times K_r$.

Consider the optimization problem 
$$
\begin{array}{rl}
\min_x  & f(\x) \\
\mbox{s.t.} &  \g(\x)=\bz, \\
& \x\in K, 
\end{array}
$$
where
$K$ is a closed convex cone that has a barrier function $\phi(\x)$.
This problem may be solved as follows.  Let $\mu_1\equiv\mu_{\rm init},\mu_2,\ldots$ be a decreasing
sequence of positive parameters tending to $0$.  Then for $k=1,2\ldots$ we solve
\begin{equation}
\begin{array}{rl}
\min_\x & f(\x)+\mu_k\phi(\x) \\
\mbox{s.t.} & \g(\x)=\bz,
\end{array}
\label{eq:primbar0}
\end{equation}
an equality-constrained optimization problem, and we define $\x_k$ to be the
optimizer or approximate optimizer.  On iteration $k+1$, we use $\x_k$
as the initial guess for the optimization algorithm, which is commonly
Newton's method.
In other words, we iteratively
solve a sequence of equality-constrained optimization problems.  
This method is called
a {\em primal} or {\em primal-only} interior-point method.  

The case most
commonly studied in the literature is the case when $f(\x)$ is
the linear function $\c^T\x$ and 
convex and the equality constraints are linear: $\g(\x)\equiv A\x-\b$ for
some matrix $A$ and vector $\b$.  In this case, there is an extensive
theory guaranteeing convergence to a global optimizer for the above
algorithm for a suitable sequence of weights $\mu_1,\mu_2,\ldots$.
See, e.g., \cite{NemNest}.

We can also define a {\em primal-dual} method as follows.
We first write the first-order (Lagrange or KKT)
optimality condition for \eqref{eq:primbar0},
which is,
\begin{equation}
\nabla f(\x^*)+\mu_k\nabla \phi(\x^*)+J(\x^*)^T\blambda=\bz,
\label{eq:opt1}
\end{equation}
where $\x^*$ is the optimizer of \eqref{eq:primbar0}, $J(\x^*)$ denotes
the first derivative (Jacobian matrix) of $\g(\x)$, and $\blambda$ is
the Lagrange multiplier.   

Assume $k$ in the previous item is fixed for now.
In the
case of $K=\R_+^n$ so that $\phi_{\rm NNO}(\x)=-\sum_{i=1}^n\log(x_i)$
and 
$$\nabla \phi(\x)=\left(
\begin{array}{c}
-1/x_1 \\
\vdots\\
-1/x_n
\end{array}
\right),$$
we define dual variable $s_i=\mu_k/x_i$ for $i=1,\ldots,n$.  
Then the optimality condition \eqref{eq:opt1},
combined with feasibility and with a rearrangement of the definition of
$s_i$ yields the following system of equations:
\begin{eqnarray*}
\nabla f(\x) - \s + J(\x)^T\blambda &=&\bz, \\
\g(\x) & = & \bz, \\
x_is_i & = &\mu_k,\quad \forall i=1,\ldots,n.
\end{eqnarray*}
The final group of equations is called ``complementarity''.  It may
be rewritten $\x\circ\s=\mu_k\e$, where ``$\circ$'' is called the {\em Jordan product}
for $\R^n_+$, and
$\e$ denotes the vector of all 1's.
The Jordan product is defined exactly as: the $i$th entry
of $\x\circ\s$ is $x_is_i$.  This notation also implies that $\x,\s$ are in the
cone, i.e., $x_i\ge 0$ and $s_i\ge 0$ for all $i=1,\ldots,n$.

In the case of $C_2^{n}$, the gradient of the barrier function is
$$\nabla \phi_{\rm SOC}(\x)=\left(
\begin{array}{c}
-x_1/d \\
x_2/d \\
\vdots\\
x_n/d
\end{array}
\right),$$
where $d=x_1^2-x_2^2-\cdots-x_n^2$.  Then we define dual variables
$s_1=\mu_k x_1/d$, $s_2=-\mu_kx_2/d$, \ldots, $s_n=-\mu_kx_n/d$.  Note that,
assuming $\x\in C_2^n$, it also follows from these formulas
that $\s\in C_2^n$.
In this
case, the optimality condition \eqref{eq:opt1} plus feasibility and the
definition of $\s$ can be written:
\begin{eqnarray*}
\nabla f(\x) - \s + J(\x)^T\blambda &=&\bz, \\
\g(\x) & = & \bz, \\
\x\circ \s & = &\mu_k\e.
\end{eqnarray*}
Here, for $C_2^n$, the Jordan product $\x\circ \s$ is defined by
$$(\x\circ\s)_i = \left\{
\begin{array}{ll}
\x^T\s, & i=1,\\
x_1s_i+s_1x_i, & i=2,\ldots, n.
\end{array}
\right.
$$
Here, $\e=[1,0,\ldots,0]^T$.  
It can be checked that $\x\circ\s=\mu_k\e$ iff $s_1=\mu_kx_1/d$ and
$s_i=-\mu_kx_i/d$ for $i=2,\ldots,n$, where $d$ is as above, provided that
$\x,\s\in C_2^n$.

Finally, if $K_1,\dots,K_r$ are all convex cones each with a Jordan product
and with Jordan identities $(\e_1,\ldots,\e_r)$,
then the Jordan identity for $K_1\times \cdots \times K_r$
is $(\e_1,\ldots,\e_r)$ and the Jordan product
is elementwise: $(\x_1,\ldots,\x_r)\circ(\s_1,\ldots,\s_r)=(\x_1\circ\s_1,
\ldots, \x_r\circ \s_r)$.

The {\em primal-dual interior-point method} consists of solving the
system of nonlinear equations:
\begin{eqnarray*}
\nabla f(\x) - \s + J(\x)^T\blambda &=&\bz, \\
\g(\x) & = & \bz, \\
\x\circ \s & = &\mu_k\e,
\end{eqnarray*}
whose variables are $(\x,\blambda,\s)$,
using Newton's method for a sequence of decreasing $\mu_k$'s, and using
the converged (or approximately converged) solution $(\x_{k-1},\blambda_{k-1},\s_{k-1})$
as the starting guess for the $k$th iteration.  
If Newton's method is applied directly to the above system,
this
yields a step called the AHO direction
(Alizadeh-Haeberly-Overton). See \cite{Todd} for an in-depth discussion.

In the nonconvex case, the known theorems are considerably weaker.
An analysis of a primal-dual interior-point method for nonconvex second-order
cone programming is presented by Yamashita and
Yabe \cite{YamashitaYabe}.  The main innovation
in that work is a merit function that ensures convergence.  
We have experimented with a primal-dual interior point method but
have not used it herein because it sometimes failed to converge
to a solution and instead became trapped close to a boundary of the
feasible region.  The hypotheses of the Yamashita-Yabe method do
not hold for the problem herein.
Therefore, the method used in our solver is a primal-only method.
However, we use the primal-dual
formulation in the computations of energy balance detailed below
in Section~\ref{sec:energy_balance}.

\section{Finite element discretization}
\label{sec:fem}

As mentioned earlier, we
assume a physical domain $\Omega\subset \R^{n_{\rm dim}}$ ($n_{\rm dim}=2$ or 
$n_{\rm dim}=3$), which is the closure of an open, bounded set with a piecewise
smooth boundary.
In this section, we describe the notation used 
to define a finite-element discretization of $\Omega$ and
$u$.
Assume that $\Omega$
is meshed with a triangulation $\mathcal{T}$.  The triangulation
is assumed to be simplicial, although the method can be extended
to meshes with hanging nodes.  The $n_{\rm dim}$-dimensional elements of this
mesh are referred to as {\em bulk elements}.

As mentioned in the introduction, we further assume that $\Omega$
contains a union $\mathcal{S}$ of curves or surfaces to represent possible
crack paths.  For the remainder of this work,
we take 
$\mathcal{S}$ to be 
the union of nonexterior
bounding curves or surfaces of the bulk
elements.
The cohesive method inserts {\em interface elements} along triangle
edges ($n_{\rm dim}=2$) or facets ($n_{\rm dim}=3$) for every pair of adjacent
bulk elements.  Let the size of $\mathcal{S}$ (number of curves or
facets) be denoted $n_{\rm e}$.  Each bulk element has its own nodes, i.e., no
node belongs to more than one bulk element.
Two adjacent bulk elements
$t_1,t_2\in\mathcal{T}$ that border on the same interface element $e$ each have nodes
in common with $e$.  Therefore, the connectivity of the mesh is determined by
nodes shared between bulk and interface elements.
Let $n_0$ denote the total number of nodes of bulk elements.
Let $n_x$ denote the number of nodal degrees of freedom
not constrained by displacement or velocity boundary conditions.
Thus, $n_x\le n_{\rm dim}n_0$.

On each time step, an optimization problem is solved
to determine the displacements at the midpoint of the time
interval. For the rest of this discussion, assume the
time step is fixed so that we omit the subscript
for time. 
Let $\u\in\R^{n_{\rm dim}\cdot n_0}$ be the vector of all nodal displacements
($n_{\rm dim}$ coordinate entries for each of $n_0$ nodes). 
As in \cite{optiIJF}, this
value of $\u$ plays the role of the unknown at the midpoint of a time-step.

Let $\x\in\R^{n_x}$ reparameterize $\u$: $\x$
stands for the degrees of freedom associated with unconstrained
nodal displacements.  The relationship between $\u$ and $\x$ is
as follows.
There is a fixed $n_{\rm dim}n_0\times n_x$ matrix $R$ such that $\u=R\x+\u_{BC}$.
Here, $\u_{BC}$ is the
$n_{\rm dim}n_0$-vector that carries information about 
displacement boundary conditions.
Note
that $\u_{BC}$ will depend on the time-step index in the
case of velocity boundary conditions.

The strain energy associated with a bulk element $t\in\mathcal{T}$ is
given by an elastic or hyperelastic energy functional.
For example, in the $n_{\rm dim}=2$ case, one choice
for the energy is the one proposed by
Knowles and Sternberg 
\cite{KnowSter1983} for
plane stress given by
\begin{equation}
  \Psi(u) = c_1\left[\trace(C) +\mathcal{J}^{-2\beta}(1+1/\beta)\right]
  \label{eq:knowles}
\end{equation}
where $\trace$ stands for the trace operator, 
$\mathcal{J}=\det(C)^{(1+1/\beta)/2}$, $C=F^TF$ (Cauchy-Green
strain), $F$ is the (2-dimensional)
displacement gradient, $c_1,\beta$ are material constants.  Inelastic bulk
material behavior is not considered herein.

The strain energy in the bulk is discretized as a function $b_0(\u)$
using quadrature over elements of $\mathcal{T}$, i.e.,
$$b_0(\u)=\int_\Omega \Psi(u_h)\,dV,$$
where $u_h$ is the finite-element interpolant specified by nodal values in the
vector $\u$.
This function $b_0(\u)$
is rewritten as $b(\x)$ 
(i.e., $b(\x)\equiv b_0(R\x+\u_{BC})$).

The momentum energy term $m_0(\u)$ arising from the
implicit midpoint rule
is derived in
\cite{optiIJF} to
be 
$$m_0(\u)= \frac{2}{\Delta t^2}
(\u - \u^i -\v^i\Delta t/2)^TM(\u - \u^i -\v^i\Delta t/2),$$
where $\Delta t$ is the time
step, $M$ is the $n_{\rm dim}n_0\times n_{\rm dim}n_0$ positive definite mass matrix,
and $\u^i$ and $\v^i$ are displacement and velocity vectors from the
preceding time step.
Define $m(\x)=m_0(R\x+\u_{BC})$.  This may be loosely
regarded as the discretization of kinetic energy;
see
\cite{optiIJF} for a more precise explanation.
Note that $m(\x)$ is a convex
quadratic function of $\x$.

Next, define an interface potential to stand for the second term of
\eqref{eq:twoterm}
as in \cite{optiIJF}.
This potential for a given element edge/surface $e\in\mathcal{S}$ is
\begin{equation}
\int_{\eta\in \Delta} g(\delta(\theta_e(\eta));d(\theta_e(\eta)))\theta_e'(\eta)\,d\eta
\label{eq:integint}
\end{equation}
where $\theta_e$ parameterizes the edge/surface $e$ with parameter
$\eta$ (a scalar for $n_{\rm dim}=2$; a 2-vector for $n_{\rm dim}=3$)
which lies in a reference domain $\Delta$,
$\delta(\cdot)$ is the
{\em effective opening
displacement} as calculated from  displacement jump in the
element boundaries, $g$ is the interface
energy function and $d$ is a damage variable
discussed below.  We follow the commonplace definition
similar to Ortiz and Pandolfi \cite{OrtizPan}:
\begin{equation}
\delta(x)=\sqrt{\llbracket u_n(x)\rrbracket^2+
(\beta^{\rm MIX})^2\Vert \llbracket u_s(x)\rrbracket\Vert^2},
\label{eq:deltaexi}
\end{equation}
where 
$u_n(\cdot)$ and $u_s(\cdot)$ are the normal and tangential
opening
displacements at a point $x\in \Omega$ that
lies on an interface.
We return to these functions below. Here, $\beta^{\rm MIX}$ is
a material constant called the {\em mixity parameter}.  Thus,
the second term $\Phi(\llbracket u\rrbracket)$ of
\eqref{eq:twoterm} is the composition
of the function $g(\delta;d)$ appearing in 
\eqref{eq:integint} with the function $\delta(\llbracket u\rrbracket)$ 
appearing in \eqref{eq:deltaexi}.  The dependence on $d$ is discussed
below.

The simplest physically reasonable choice for $g$ prior to the
introduction of damage is:
\begin{equation}
  g(\delta) = \left\{
  \begin{array}{ll}
    l\delta+q\delta^2, & \delta\in[0,\delta_u] \\
    l\delta_u + q\delta_u^2, & \delta \ge \delta_u,
  \end{array}
  \right.
\end{equation}
where $\delta_u$, the ultimate opening displacement, is a material parameter; 
quadratic coefficient 
$q=-\sigma_c/(2\delta_u)$,
$l=\sigma_c$, and $\sigma_c$, the critical traction, is another
material parameter.  One checks that with these formulas, $g$ is a piecewise
$C^1$ (continuous function and first derivative) quadratic function.
Its first derivative with respect to $\delta$, $g'$, is therefore piecewise linear and continuous.
With these
choices of the three parameters, $g'(0)=\sigma_c$, indicating that 
the initial traction (first derivative of energy with respect to $\delta$)
is $\sigma_c$.
The area under the curve of the plot of $g'(\delta)$
is $G_c=\sigma_c\delta_u/2$, a material parameter called the ``critical energy
release rate''.
(Note: of the three material parameters $\sigma_c,\delta_u,G_c$, only two
can be chosen independently as the previous equality demonstrates.)

\begin{figure}[htp]
\begin{center}
\epsfig{file=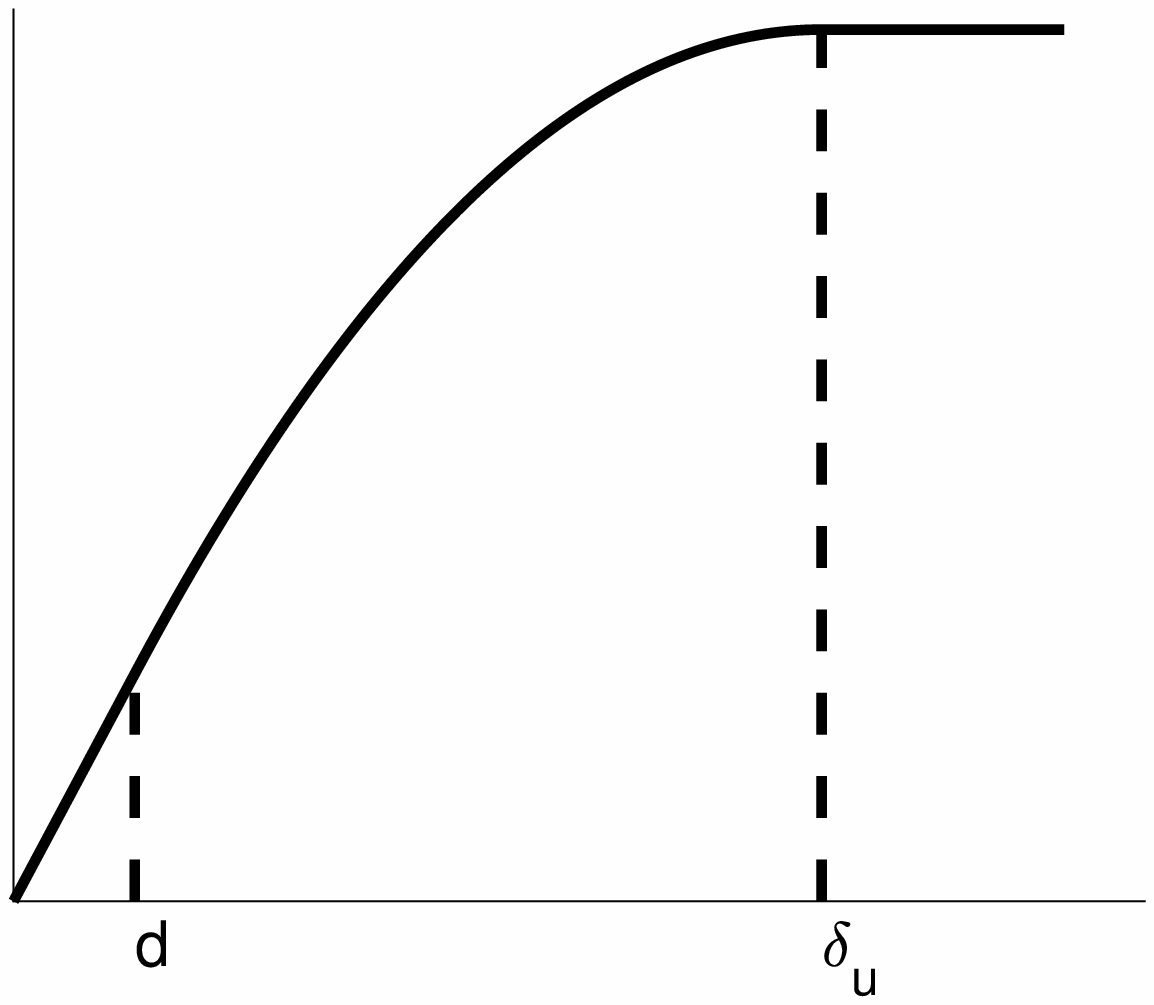,height=5cm} 
\epsfig{file=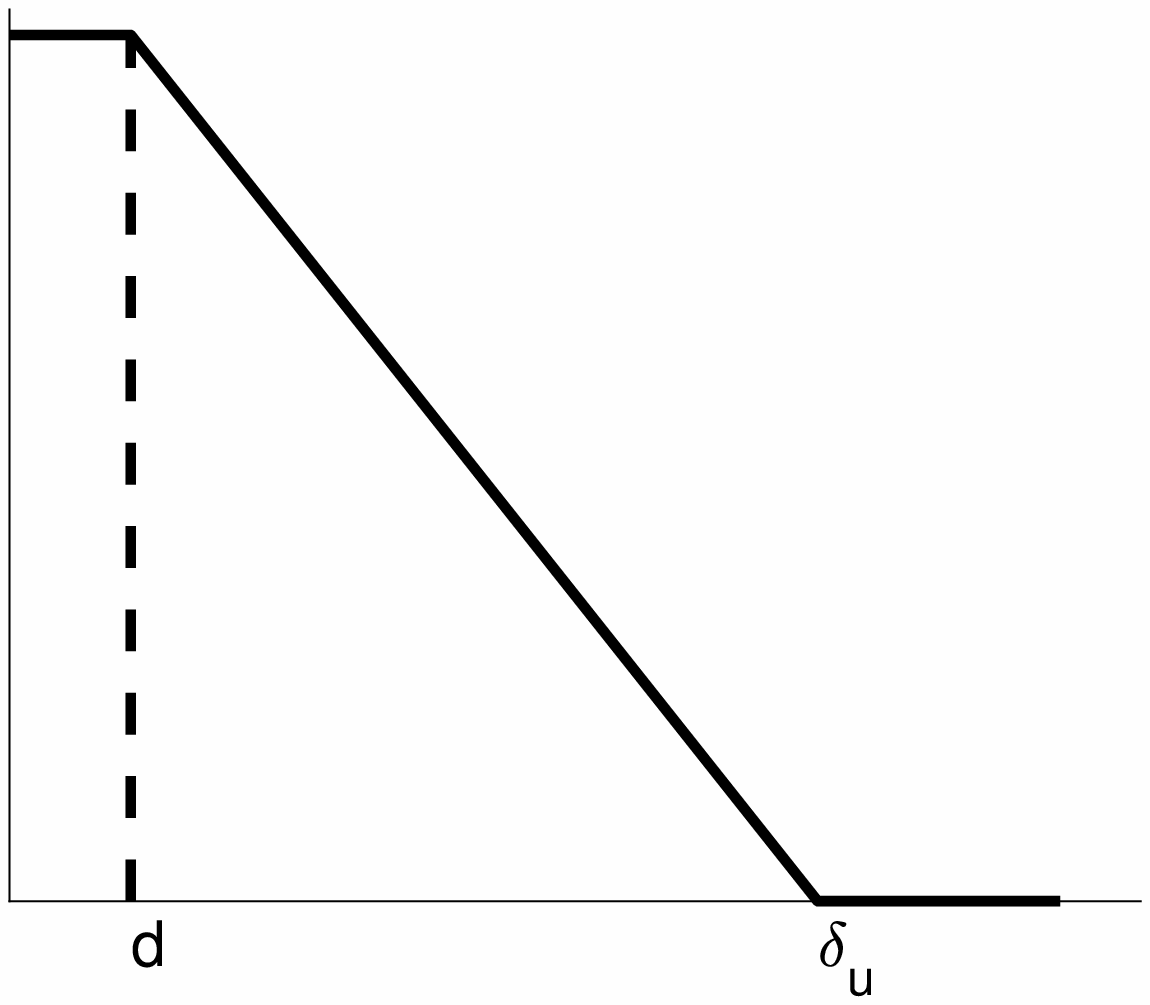,height=5cm} 
\end{center}
\caption{Plot of the function $g(\delta;d)$ (left) 
and $g'(\delta;d)$ (right) defined
by \eqref{eq:genergy}.  Although difficult to discern at this scale, the segment
of the curve on the left corresponding to abscissas lying in $[0,d]$ is straight (linear).}
\label{fig:gplot}
\end{figure}

We now extend this formula to include a nonnegative
scalar damage parameter $d$,
initially equal to 0.
The role of scalar $d$ is to model irreversible
damage to the interface.  As in \cite{OrtizPan} and many other
previous works, we define this parameter
equal to the maximum opening displacement (not
exceeding $\delta_u$)
encountered over previous time values.   When $d=\delta_u$, the
interface has no remaining cohesion.  The extended formula
is:
\begin{equation}
g(\delta;d)=\left\{
\begin{array}{ll}
l(d)\delta, & \delta\in [0,d], \\
l(d)\delta + q(\delta-d)^2 & \delta \in [d,\delta_u], \\
l(d)\delta_u +q(\delta_u-d)^2 & \delta >\delta_u, 
\end{array}
\right.
\label{eq:genergy}
\end{equation}
where now $q=-\sigma_c/(2\delta_u)$,
$l(d)=-2(\delta_u-d)q$.  When $d=0$, the formula in the previous
paragraph is recovered.
See Fig.~\ref{fig:gplot}.  Unlike previous works such as \cite{OrtizPan},
this formula implies
that the material retains a residual critical stress when $0<d<\delta_u$.
In contrast, most previous works specify that after the onset of damage,
the interface behaves like an initially elastic interface (i.e., unloads
to the origin). The difference
in practice between the two models appears to be minor, but our formulation
has the mathematical advantage
that it prevents a pathological situation in which the second
derivative of $g$ can have an unboundedly large value.

In the finite element
approximation, the integral \eqref{eq:integint} is computed for each edge ($n_{\rm dim}=2$) or
facet ($n_{\rm dim}=3$) 
with a $n_{\rm g}$-point Gauss quadrature rule.
Let $n_{\rm i}=n_{\rm e}n_{\rm g}$ stand for the total number
of Gauss-points of interfaces.
Let us introduce a new variable
$\s_0\in\R^{n_{\rm i}}$ that represents the effective opening displacements
at the Gauss points.  In other words, $\s_0$ stores the vector of values of
$\delta(\theta_e(\eta_\iota))$, $\iota=1,\ldots,n_{\rm g}$,
described above at 
each Gauss point 
$\eta_1,\ldots,\eta_{n_{\rm g}}$
of each interface $e$.
Then the potential due to interfaces is written
\begin{equation}
h(\s_0;\d)
=\sum_e\sum_{\iota=1}^{n_{\rm g}}\omega_{e,\iota}g((s_0)_{e,\iota};d_{e,\iota});
\label{eq:hs0}
\end{equation}
this is the finite-element approximation to $\Phi(\llbracket u\rrbracket )$ that appears in \eqref{eq:twoterm}.
We have associated a damage variable
$d_{e,\iota}$ for $e\in \mathcal{S}$ and $\iota=1,\ldots,n_{\rm g}$, 
i.e., one for each of the $n_{\rm i}$ interface Gauss points, that
indicate the level of damage.
Here, $\omega_{e,\iota}$ is the quadrature weight.  This function $h$
is separable in the entries of $\s_0$, and hence its Hessian is a diagonal matrix.

Let $\s_1\in\R^{n_{\rm i}}$ denote the
normal opening displacement. Let $\s_2,\ldots,\s_{n_{\rm dim}}\in\R^{n_{\rm i}}$ 
stand for the tangential
opening displacements scaled by $\beta^{\rm MIX}$.  (Note that for three-dimensional
problems, there is not a unique way to define a tangential coordinate
system at each point on an interface.  Any method for defining tangential
coordinates is acceptable provided that it is applied consistently for
the duration of the solution procedure.)
Next, we revisit \eqref{eq:deltaexi}.  First, the left-hand side is needed
only at the Gauss points of the interfaces.  In this case, the left-hand side is
just $(s_0)_{e,\iota}$.  A similar substitution may be made on the right-hand
side, so this constraint is rewritten as
$$(s_0)_{e,\iota}=\sqrt{(s_1)_{e,\iota}^2+\cdots+(s_{n_{\rm dim}})_{e,\iota}^2},$$
for $e=1,\ldots,n_{\rm e}$, $\iota=1,\ldots,n_{\rm g}$.  Observe that
this constraint is nondifferentiable at the origin, that is,
when for some $e,\iota$, 
$(s_1)_{e,\iota}=\ldots=(s_{n_{\rm dim}})_{e,\iota}=0$.  This nondifferentiability
is fundamental to the model and is precisely the reason
why it is able to capture the initially rigid interface behavior.
A detailed explanation of the role of this nondifferentiability
is provided in \cite{optiIJF}.

The vectors $\s_1,\ldots,\s_{n_{\rm dim}}$ containing the
components of the opening displacements are functions of the 
displacements stored in $\u$.  In other words, we can
determine entries of $\s_1,\ldots,\s_{n_{\rm dim}}$ by evaluating
jumps of displacements interpolated from shape functions.
As mentioned earlier, we have reparameterized $\u$ by $\x$.
Therefore we can define geometric
functions $c_{k,e,\iota}$ such that
\begin{equation}
(s_k)_{e,\iota}=c_{k,e,\iota}(\x),
\label{eq:sxconstr}
\end{equation}
for  $k=1,\ldots,n_{\rm dim}$, 
$e=1,\ldots,n_{\rm e}$, $\iota=1,\ldots,n_{\rm g}$.  
As observed in \cite{optiIJF},
these functions are nonlinear
because of geometric nonlinearity, namely, the
normal and tangent directions depend on the current values of 
the displacements.  
In addition to geometry, the functions
$c_{2,e,\iota},\ldots,c_{n_{\rm dim},e,\iota}$ also have the mixity factor
$\beta^{\rm MIX}$ encoded in them.

As mentioned in the introduction, one advantage of the
interior-point formulation is the ability to handle
conic convex constraints essentially for free. 
One special case of conic convex
constraints is linear inequality contraints.  
Let us assume that the system has additional
linear constraints of the form $E\x\ge \a$.
(Later
on, we will use these inequalities to model a simple form 
of a contact constraint.)
Here, $E\in\R^{n_{\rm LI}\times n_x}$ is a known
matrix and $\a\in\R^{n_{\rm LI}}$ is a known vector, both of which
may vary from one time-step
to the next, and $n_{\rm LI}$ denotes the number of 
linear inequality constraints.

Thus, the optimization problem to solve for one time-step in 
the model of cohesive fracture is:
\begin{equation}
\begin{array}{rll}
\min_{\x,\s_0,\s_1,\ldots,\s_{n_{\rm dim}}} & m(\x)+b(\x)+h(\s_0;\d)+\f^T\x \\
\mbox{s.t.} 
& (s_0)_{e,\iota}=\sqrt{(s_1)_{e,\iota}^2+\cdots+(s_{n_{\rm dim}})_{e,\iota}^2}& \forall e,\iota, \\
& (s_k)_{e,\iota}=c_{k,e,\iota}(\x) & \forall k=1,\ldots,n_{\rm dim},\forall e,\forall\iota, \\
& E\x\ge \a,\\
& (s_1)_{e,\iota}\ge 0 & \forall e,\forall\iota. \\
\end{array}
\label{eq:optmodeldisc}
\end{equation}
The terms in the objective have already been discussed except for
the last term, which stands for the sum of traction and body forces
from \eqref{eq:twoterm}.
The first, second and third constraints were
already discussed.  The fourth
prevents interpenetration between neighboring
elements.

\section{Interior-point method for fracture}
\label{sec:ipfm}

A key modification to
\eqref{eq:optmodeldisc}
that makes it amenable to 
an interior-point method is to replace the first equality constraint with 
an inequality constraint:
\begin{equation}
\begin{array}{rll}
\min_{\x,\s} & m(\x)+b(\x)+h(\s_0;\d)+\f^T\x \\
\mbox{s.t.}
& (s_0)_{e,\iota}\ge\sqrt{(s_1)_{e,\iota}^2+\cdots+(s_{n_{\rm dim}})_{e,\iota}^2}& \forall e,\forall\iota, \\
& (s_k)_{e,\iota}=c_{k,e,\iota}(\x) & \forall k=1,\ldots,n_{\rm dim},\forall e,\forall\iota, \\
& E\x\ge \a, \\
& (s_1)_{e,\iota}\ge 0 & \forall e,\iota.
\end{array}
\label{eq:mainprob}
\end{equation}
Replacing the equality by an inequality constraint
does not change the optimizer because $h$ is a nondecreasing
function of $\s_0$. (In fact, it is perturbed to
a strictly increasing function as described below.)
This implies that the optimal solution to \eqref{eq:mainprob}
satisfies
$(s_0)_{e,\iota}=\sqrt{(s_1)_{e,\iota}^2+\cdots+(s_{n_{\rm dim}})_{e,\iota}^2}$.

The benefit of this change is that an
equality constraint
of the form $(s_0)_{e,\iota}=\Vert(\s_{1:n_{\rm dim}})_{e,\iota}\Vert$ defines a nonconvex set
with a complicated (nonmanifold) structure, whereas 
the constraint $(s_0)_{e,\iota}\ge\Vert(\s_{1:n_{\rm dim}})_{e,\iota}\Vert$
defines a convex set $C_2^{n_{\rm dim}+1}$, the
second-order cone.  Here, $(\s_{1:n_{\rm dim}})_{e,\iota}$ is the vector in 
$\R^{n_{\rm dim}}$ whose
entries are $((s_1)_{e,\iota},\ldots,(s_{n_{\rm dim}})_{e,\iota})$.

As discussed in Section~\ref{sec:ipbg},
the primal-only interior-point method replaces the inequality
constraints in the preceding formulation with log-barrier terms.  The
parameter $\mu>0$ starts at a large value and decreases to close
to zero.  This leads to the following formulation:
\begin{equation}
\begin{array}{rl}
  \min_{\x,\s} & m(\x)+b(\x)+h_{\mu}(\s_0;\d)+\f^T\x
  +\mu\phi_{\rm NNO}(E\x-\a)
  \\
& \quad \mbox{} +
\mu\sum_{e=1}^{n_{\rm e}}\sum_{\iota=1}^{n_{\rm g}}
\zeta_{e,\iota}\left[\phi_{\rm SOC}((\s_{0:n_{\rm dim}})_{e,\iota}) +\phi_{\rm NNO}((s_1)_{e,\iota})\right] \\
\mbox{s.t.}
& (s_k)_{e,\iota}=c_{k,e,\iota}(\x) \quad \forall k=1,\ldots,n_{\rm dim},\forall e,\forall\iota. \\
\end{array}
\label{eq:primbar}
\end{equation}
The positive weight $\zeta_{e,\iota}$ is
defined by \eqref{eq:zetadef} in the next section.
In addition, we have perturbed $h$ of
\eqref{eq:mainprob} to
$h_\mu$ in \eqref{eq:primbar}, which is defined by \eqref{eq:hmu} 
in the next section.  As $\mu\rightarrow 0$, $h_\mu\rightarrow h$,
thus recovering the original problem.

This formulation still contains equality constraints, but they are easily eliminated
via substitution.  In particular,
we substitute $c_{k,e,\iota}(\x)$ for $(s_k)_{e,\iota}$ 
(as in \eqref{eq:sxconstr})
thus
eliminating $(s_k)_{e,\iota}$, $k=1,\ldots,n_{\rm dim}$, 
$e=1,\ldots,n_{\rm e}$, $\iota=1,\ldots,n_{\rm g}$.  This elimination leaves
the following unconstrained problem:
\begin{equation}
  \begin{array}{rl}
\min_{\x,\s_0} & m(\x)+b(\x)+h_{\mu}(\s_0;\d)+\f^T\x
+\mu\phi_{\rm NNO}(E\x-\a) \\
 & \quad \mbox{} +
\mu\sum_{e=1}^{n_{\rm e}}\sum_{\iota=1}^{n_{\rm g}}
\zeta_{e,\iota}\bigg[\phi_{\rm SOC}((s_0)_{e,\iota},\c_{1:n_{\rm dim},e,\iota}(\x)) \\
&\quad \quad \mbox{}
  +\phi_{\rm NNO}(c_{1,e,\iota}(\x))\bigg]
  \end{array}
  \label{eq:primbar2}
\end{equation}

The interior-point code solves \eqref{eq:primbar2}
via a trust-region method.  This is a standard technique
to extend Newton's method to 
nonconvex (unconstrained) optimization.  More details on the method are
in Section~\ref{sec:compproc}.

\section{Sources of nonconvexity}
\label{sec:nonconv}
The optimization problem \eqref{eq:mainprob} contains three
sources of nonconvexity as follows. The bulk energy $b(\x)$  is
nonconvex in the displacements for nonlinear hyperelasticity,
$h(\s_0;\d)$ is nonconvex in $\s_0$ due
to the negative coefficient of the quadratic term in
\eqref{eq:genergy}, and
the constraints
$(s_i)_{e,\iota}=c_{i,e,\iota}(\x)$, $i=1,\ldots,n_{\rm dim}$,
are nonlinear.
(Note that in passing from \eqref{eq:optmodeldisc} to
\eqref{eq:mainprob}, a fourth source of nonconvexity was
eliminated by replacing an equality with an inequality.)
Of these three sources,
the nonconvexity of $h$ is the most challenging
to handle, and
it is also the most fundamental to the application.  Convexity
of $b(\x)$ could be recovered by adopting a simpler mechanical
model such as linear elasticity. Linearity in
the constraint $(s_k)_{e,\iota}=c_{k,e,\iota}(\x)$
could be recovered by simply assuming that the normal
vectors to the interfaces are determined by the initial rather
than current configuration.  
There is, however, no apparent
way to replace or approximate $h$ with a convex function
because $g$ must have a nonconvex form similar to the
form depicted in
Fig.~\ref{fig:gplot}(left) to be physically meaningful.

Interior-point methods for nonconvex problems are considerably
more delicate than for convex problems, and we were
required to implement several stabilization methods in the
interior-point framework to cope with the nonconvexity of $h$,
which are as follows.

\begin{enumerate}
\item
The interior-point method converges significantly faster if we weight
the log-barrier terms for the
constraints associated with interfaces by a factor proportional to
their ``local'' length.  In particular, we introduced weight
$\zeta_{e,\iota}$ in \eqref{eq:primbar}, which is defined as
\begin{equation}
  \zeta_{e,\iota}=10^4G_c\omega_{e,\iota},
  \label{eq:zetadef}
\end{equation}
where $\omega_{e,\iota}$ is the quadrature weight in \eqref{eq:hs0}.
Using a weight proportional to $\omega_{e,\iota}$ 
is physically natural because it means, for example, that
the contribution to the barrier function from an interface is invariant
(up to discretization error) if the interface is subdivided into
smaller pieces.  Making the weight proportional to $G_c$, the critical
energy release rate, is also natural since the other terms in the objective
function stand for work or energy quantities.
Weighting is not necessary (and is typically not even
considered) in the case of convex interior-point methods.

\item
On intermediate stages of the interior-point method, the interfaces
are favored to open by the log-barrier term
associated with the constraint
$(s_1)_{e,\iota}\ge 0$ when $\mu$ is large.  
Without extra measures, they can open
by more than $\delta_u$, in which case their traction is 0 and they no longer
hold the body together.  Their traction is not recovered as $\mu$
is decreased due to the nonconvexity of $h$.
For this reason, a quadratic regularization term is added to the cohesive
traction.  This term has the form 
$$(5\cdot 10^5)\alpha\omega_{e,\iota}\max(1-d_{e,\iota}/\delta_u,8\cdot 10^{-6})(s_0)_{e,\iota}^2,$$
where $\alpha>0$ is a small scalar, $\omega_{e,\iota}$ is the
quadrature weight in \eqref{eq:hs0},
$d_{e,\iota}$ is the damage value, $(s_0)_{e,\iota}$ is the (unknown)
effective opening displacement of 
Gauss point $\iota$ of the interface $e$.
From now on, we denote:
\begin{equation}
h_\alpha(\s_0;\d)=h(\s_0;\d) + \sum_{e,\iota} (5\cdot 10^5)\alpha\omega_{e,\iota}
\max(1-d_{e,\iota}/\delta_u,0.000008)(s_0)_{e,\iota}^2, 
\label{eq:hmu}
\end{equation}
As already noted in
\eqref{eq:primbar2},
the method sets $\alpha$ in \eqref{eq:hmu} equal to the barrier parameter
$\mu$ (so that the original problem is
recovered as $\mu\rightarrow 0$).

\item
  The Newton step associated with the interior-point step is not always
  well defined  because the Hessian may not be positive definite.  To address
  this problem, we use a common regularization to extend
  Newton to nonconvex optimization, namely, the trust-region
  method.  Some modifications to the standard trust-region method were necessary
  for this application; these are described in detail in Section~\ref{sec:compproc}.

\end{enumerate}
The net effect of these stabilization techniques is that the method
is reasonably fast and robust on all the computational experiments tried so far.

\section{Attaining feasibility with ``Phase I'' initialization}
\label{sec:phaseI}

An issue to address is that interior-point methods for optimization
require an interior starting point.
There are so-called ``infeasible'' 
interior-point methods that ease this restriction, 
but the theory for such methods in the case of nonconvex
problems is not well developed, and in practice they can be difficult
to use.  
A natural choice for initialization is
the converged solution for the previous time
step, but such a solution may violate the inequality constraints
such as contact.  In addition, for moving displacement boundary conditions
(i.e., velocity boundary conditions), the boundary nodes must be
displaced each iteration to new positions to attain feasibility.
Depending on the magnitude of the velocity, this can cause elements
along the boundary to become inverted or nearly inverted, which makes
the bulk material model behave poorly to the extent that recovering
a noninverted shape is unattainable with Newton's method.

One solution to the problem described in the last paragraph
is to take small load steps, thus limiting the degree of misshapenness
among boundary elements.  The drawback of this technique is that it makes
the time-step dependent on the mesh size because a finer mesh implies
that a smaller boundary distortion can be tolerated.  Recall that
the method proposed herein is an implicit method.  A big advantage
of implicit methods is exactly that the time step and mesh size can
be chosen independently.  Thus, coping with moving boundaries by subdividing
time-steps undermines one of the main benefits
of implicit methods.

An alternative approach, adopted here, is to start with a so-called ``Phase I''
initialization.  During Phase I, the optimization problem is modified 
with a new variable $t$ constrained $t\ge 0$ and a ``big $M$'' term in the objective.
A preliminary version of the modified problem is as follows. 
\begin{equation}
  \begin{array}{rl}
    \min_{\x,\s_0,t} & m(\x)+b(\x)+h_{\mu_{\rm init}\sqrt{M}}(\s_0;\d)+\f^T\x +Mt \\
    &\quad\mbox{}
+\mu_{\rm init}\phi_{\rm NNO}(E\x+t\e-\a)+\mu_{\rm init}\phi_{\rm NNO}(t) \\
 & \quad \mbox{} +
\mu_{\rm init}\sum_{e=1}^{n_{\rm e}}\sum_{\iota=1}^{n_{\rm g}}
\zeta_{e,\iota}\bigg[\phi_{\rm SOC}((s_0)_{e,\iota},\c_{1:n_{\rm dim},e,\iota}(\x)) \\
&\quad \quad \mbox{}
  +\phi_{\rm NNO}(c_{1,e,\iota}(\x)+t\e)\bigg]
  \end{array}
  \label{eq:primbar2phaseI0}
\end{equation}
Here, $\e$ is the vector of all 1's.  Note that during Phase I, $\mu$ is held
at its initial value as suggested by the notation $\mu_{\rm init}$ appearing in
\eqref{eq:primbar2phaseI0}.

The rationale for \eqref{eq:primbar2phaseI0} is as follows.  The variable
$t$ can be initialized to a large positive value, which ensures that
the barrier functions $\phi_{\rm NNO}(\cdot)$ have positive arguments in all occurrences.
In addition, the variables  $(s_0)_{e,\iota}$ can be
initialized to sufficiently large values for all interface
Gauss points $(e,\iota)$ to ensure that the argument of $\phi_{\rm SOC}$ is
feasible.  Thus, except for the issue of element inversion that is discussed
below, it is straightforward to find an initial feasible point for
\eqref{eq:primbar2phaseI0}.  On the other hand, at the optimizer
for \eqref{eq:primbar2phaseI0}, $t>0$ will be a relatively small value
because the term $Mt$ in the objective function penalizes a large $t$.
If a feasible
solution to \eqref{eq:primbar2phaseI0} has a small value of $|t|$, 
then we can simply set $t=0$ and expect to have a feasible
solution for    \eqref{eq:primbar2}.

If the solution so obtained
is feasible for \eqref{eq:primbar2}, then we
proceed to a solver for
\eqref{eq:primbar2} using this initial guess.
If not, we re-solve \eqref{eq:primbar2phaseI0}
with $M$ increased by a factor of 8 (initially, $M=64$).  A larger value of $M$ means a smaller
value of $t>0$ at the optimizer.  The parameter $\alpha$
appearing in $h_\alpha$ in \eqref{eq:hmu}
is $\mu_{\rm init}\sqrt{M}$ in order to ensure that the quadratic penalty term in $h$
is not swamped by the $Mt$ term in the objective.

The ``big-$M$'' technique thus handles feasibility of the
barrier functions.  Next, we further modify \eqref{eq:primbar2phaseI0} to prevent
element inversion that can be caused by
moving boundary conditions.
Recall that the vector of all degrees of freedom (DOFs)
$\u$ is parameterized by $\u=R\x+\u_{\rm BC}$, where $\u_{\rm BC}$ carries
the prescribed displacements and $R$ is a fixed matrix.  There is always
a mathematically equivalent way to
to write the boundary constraints as $B\u=\b$, i.e., linear equations that must be satisfied
by $\u$.  The transformation from $(R,\u_{\rm BC})$ to
$(B,\b)$ can be carried out using standard numerical linear algebra
such as QR factorization.   Note that both $\u_{\rm BC}$ and $\b$ may depend
on the time step $\tau$, whereas $R$ and $B$ are fixed throughout the computation.

In Phase I, we let the vector of unknown displacements be $\u$ rather
than $\x$. Instead of equality constraints $B\u=\b$, we introduce
two inequalities $B\u+t\e\ge\b$
and $B\u-t\e\le \b$, where $t$ is the Phase-I artificial
variable.  In turn, these two constraints are
replaced by additional terms in the barrier objective function
of the form $\mu\phi_{NNO}(B\u+t\e-\b)$ and $\mu\phi_{NNO}(-B\u+t\e+\b)$.
Thus, the problem is reformulated as:
\begin{equation}
  \begin{array}{rl}
    \min_{\u,\s_0,t} & m_0(\u)+b_0(\u)+h_{\mu_{\rm init}\sqrt{M}}(\s_0;\d)+\f_0^T\u +Mt \\
    &\quad\mbox{}
+\mu_{\rm init}\phi_{\rm NNO}(E_0\u+t\e-\a_0)+\mu_{\rm init}\phi_{\rm NNO}(t) \\
 & \quad \mbox{} + \mu_{\rm init}\phi_{\rm NNO}(B\u+t\e-\b)+\mu_{\rm init}\phi_{\rm NNO}(-B\u+t\e+\b) \\
 & \quad \mbox{} +
\mu_{\rm init}\sum_{e=1}^{n_{\rm e}}\sum_{\iota=1}^{n_{\rm g}}
\zeta_{e,\iota}\bigg[\phi_{\rm SOC}((s_0)_{e,\iota},\c_{0;1:n_{\rm dim},e,\iota}(\u)) \\
&\quad \quad \mbox{}
  +\phi_{\rm NNO}(c_{0;1,e,\iota}(\u)+t\e)\bigg]
  \end{array}
  \label{eq:primbar2phaseI}
\end{equation}
Here, we recall the notation introduced earlier that $b_0(R\x+\u_{\rm BC})\equiv b(\x)$
and similarly for $m_0(\cdot), \f_0, E_0, \a_0, \c_0(\cdot)$.  Now it is apparent
that by selecting both $t$ and $(\s_0)_{e,\iota}$ sufficiently large, all barrier
constraints are feasible and no element is inverted by the initial guess for $\u$
(the solution from the previous time-step).
It is also apparent that if $t$ is driven to a sufficiently small positive number
by optimization on  \eqref{eq:primbar2phaseI}, then a feasible solution for the
original (phase II) problem is obtained.

\section{Computational procedure}
\label{sec:compproc}
In this section, we provide further details of the computational procedure
of the cohesive solver.  As mentioned in the last section,
during Phase I, the displacement boundary conditions
are enforced as inequalities, and therefore the vector of unknowns
for the displacements in \eqref{eq:primbar2phaseI} is $\u$ rather than
$\x$.
Let $\Pi_\tau$ be the linear projection
that maps $\u$ to $\x$, i.e., $\Pi_\tau(\u)=\argmin_\x\Vert(R\x+\u_{BC,\tau})-\u\Vert$, 
where $\tau$ is the time step.
The feasibility test that terminates Phase I that was described
in the previous section in more detail is:
apply $\Pi_\tau$ to $\u$ that solves \eqref{eq:primbar2phaseI} to make sure that
the vector $\x$ thus obtained is feasible for the main phase.

For the purpose of notation in the solvers,
let $\bxi$ denote the concatenation $(\x,\s_0)$, the
variables of the interior method
\eqref{eq:primbar2}, or $\bar\bxi = (\u,\s_0,t)$ in the case
of Phase I when \eqref{eq:primbar2phaseI} is solved.

In the algorithms that follow,
these variables $\bxi$ and $\bar\bxi$ are not meant to stand
for new independent program variables but merely notational
shorthand.  For
example, if $\x$ is updated in a code that follows, 
then $\bxi$ is also updated implicitly
since $\bxi$ contains $\x$ as a subvector.

We make the following observation: given a value of $\x$ or $\u$,
it is possible to efficiently compute the optimal extension
of $\x$ to $\bxi$ or $\u$ to $\bar\bxi$, where ``optimal'' in this context
means minimizing the relevant objective function $f$.  One observes
from \eqref{eq:mainprob} that once $\x$  is specified, the
remaining variables $(s_0)_{e,\iota}$, $e=1,\ldots, n_{\rm e}$,
$\iota=1,\ldots,n_{\rm g}$ and $t$ (in the case of $\bar\bxi$)
are decoupled and may be optimized individually
using a univariate procedure (e.g., bisection).
Let us denote these optimal values as $\bxi^*(\x)$, $\bar\bxi^*(\u)$.

The top-level procedure is described in Fig.~\ref{fig:toplevel}.
Every third time step, the algorithm computes matrices
$H_{\mu_i}$ for $i=1,\ldots, n_\mu$, which are
positive definite matrices used in the trust-region method.

The variables maintained in the main loop
from one time step to the next are $\u$ and $\d$, which
are superscripted with the time step index $\tau$.
As discussed earlier, $\u$ (or $\x$) encodes the displacements while
$\d$ is the
damage state, which is updated from the displacements computed
on each step.

\begin{figure}[htp]
\begin{center}
\fbox{\begin{minipage}{\dimexpr\textwidth-2\fboxsep-2\fboxrule\relax}
\begin{tabbing}
  +\=++\=++\=\kill
  \> From initial conditions determine initial guess for $\u^0$. \\
  \> $\d^0:=\bz$. \\
  \> $ \{H_{\mu}\}_\mu := 0\,\forall \mu;\quad H^{\rm PHASE1}:=0.$ \\  
  \> for $\tau := 1,\ldots,n_{\rm step}$ \\
  \> \> if $\mbox{rem}(\tau,3)==1$ \\
  \> \> \> $H^{\rm PHASE1},\{H_\mu\}_\mu=\texttt{solver}(\tau,\u^{0},\d^0,H^{\rm PHASE1},\{H_{\mu}\}_\mu,
  \mbox{\sl PREPROCESS})$. \\
  \> \> end if \\
\> \> $\u^\tau :=\texttt{solver}(\tau, \u^{\tau-1},\d^{\tau-1},H^{\rm PHASE1},\{H_\mu\}_\mu,\mbox{\sl ORDINARY})$. \\
\> \> Compute $\d^\tau$ from $\d^{\tau-1},\u^{\tau}$. \\
\> end for
\end{tabbing}
\end{minipage}}
\end{center}
\caption{Top-level procedure.}
\label{fig:toplevel}
\end{figure}

The sequence of $\mu$ values used in the interior-point
methods are fixed in advance as a geometrically decreasing sequence:
$\mu_i=\mu_{\rm init}\rho_{\mu}^{i-1}$ for $i=1,2,\ldots,n_{\mu}$.
The choice of parameters used is $\mu_{\rm init}=5\cdot 10^{-5}$, $\rho_\mu=0.125$,
and $n_{\mu}=6$.  This means that the ultimate value 
is $\mu_{n_\mu}\approx 1.5\cdot 10^{-9}$.

The procedure \verb+solver+ to solve one time step
using the interior-point method is 
detailed in Fig.~\ref{fig:solvermainprob}.

\begin{figure}[htp]
\begin{center}
\fbox{\begin{minipage}{\dimexpr\textwidth-2\fboxsep-2\fboxrule\relax}
\begin{tabbing}
+\=++\=++\=++\=++\=\kill
\> FUNCTION $\texttt{solver}(\tau,\u,\d,H^{\rm PHASE1},\{H_\mu\}_\mu,\mbox{\sl flag1})$ \\
\> Determine functions $b(\cdot)$, $m(\cdot)$, $b_0(\cdot)$, $m_0(\cdot)$ for time step $\tau$. \\
\> $\bar\bxi:=\bar\bxi^*(\u)$ \\
\> $M:=64;\quad\mu:=\mu_1;\quad R:=1.0.$\\
\> {\em /* Phase 1 to find feasible $\x$ by minimizing \eqref{eq:primbar2phaseI}.}\\
\> {\em \quad Variable $M$ varies in this loop, while $\mu$ is fixed. */} \\
\> loop \\
\> \> Let $f(\cdot)$ be the objective function of \eqref{eq:primbar2phaseI}.\\
\> \> $(\u,\s_0,t, R) := \texttt{minimize}(f, \bar\bxi, H^{\rm PHASE1}, R).$ \\
\> \> $\x := \Pi_\tau(\u).$ \\
\> \> if $(\x,\s_0)$ feasible for \eqref{eq:primbar2}\\
\> \> \> break \\
\> \> end if \\
\> \> $M := 8M$. \\
\> end loop \\
\> $\bxi:=(\x,\s_0)$. \\
\> if {\sl flag1} =={\sl PREPROCESS} \\
\> \> $H^{\rm PHASE1} :=  \nabla^2 f(\u).$ \\
\> end if \\
\> {\em /* Phase 2 to minimize \eqref{eq:primbar2}.  Variable $\mu$ decreases} \\
\> {\em \quad according to fixed schedule.  */} \\
\> for each $\mu:=\mu_1,\mu_2,\ldots,\mu_{n_\mu}$ \\
\> \> Let  $f(\cdot)$ be the objective function of \eqref{eq:primbar2}. \\
\> \> if {\sl flag1} =={\sl PREPROCESS} \\
\> \>  \> Replace $h_\mu(\cdot)$ appearing in \eqref{eq:primbar2} by $h_{\mu_1}(\cdot)$.\\
\> \> end if \\
\> \> $(\bxi,R) := \texttt{minimize}(f,\bxi,H_\mu,R)$. \\
\> \> if {\sl flag1}=={\sl PREPROCESS} \\
\> \> \> $H_\mu:=\nabla^2 f(\bxi)$ \\
\> \> end if \\
\> end for \\
\>  if {\sl flag1} == {\sl PREPROCESS} \\
\> \> return $H^{\rm PHASE1}, \{H_\mu\}_\mu.$ \\
\> else \\
\> \> return $R\x+\u_{\rm BC}.$ \\
\> end if
\end{tabbing}
\end{minipage}}
\end{center}
\caption{Solver for \eqref{eq:mainprob}.}
\label{fig:solvermainprob}
\end{figure}

Note that the bulk strain-energy function $b(\x)$ in general
depends on the time step index.  This is because this function $b(\cdot)$ encodes
all of the boundary and loading conditions.  The momentum energy
term
$m(\x)$ also varies with the time step because it incorporates
values of the nodal velocities from previous time steps.

The interior-point minimizer, which appears in Fig.~\ref{fig:primalonly},
is a modification of the standard trust-region method
(see, e.g., Nocedal and Wright \cite{NocedalWright}).
As is standard for this method, the function $m$ appearing in the ratio test $\rho$ is
the quadratic model, that is,
$m(\x):=\nabla f(\bxi)^T\x + (\x^T\nabla^2f(\bxi)\x)/2$.

The two modifications are as follows.  First, instead of the usual
identity matrix added to regularize the Hessian, we use the
Hessian computed on an artificial preliminary problem plus
a multiple of the identity.  This Hessian
captures the geometry of the space of $\bxi$  more accurately than a plain identity
matrix and
therefore leads to faster convergence.  This is because the objective
function $f(\bxi)$ is highly anisotropic;
some search directions (e.g., those that create significant interpenetration)
cause a large jump in the objective, while others only a small change.

The second modification is that, in addition to the usual ratio test
on function values
for determining when to accept a step, the function also implements a ratio
test on gradient norms, namely, the variable $\rho_g$ appearing
in Fig.~\ref{fig:primalonly}.  It follows from Taylor's theorem that
the numerator of the definition of $\rho_g$ tends to zero rapidly as
$\Delta \bxi$ gets small, so a large value $\rho_g$ indicates that the
function is not behaving according to the Taylor prediction.
This modification was necessary because of
pathological cases of the
trust region method in which the objective function decreases while the gradient
blows up to infinity as the boundary of a cone is approached obliquely.
Applying a ratio test to the gradient prevents such occurrences.

\begin{figure}[htp]
\begin{center}
\fbox{\begin{minipage}{\dimexpr\textwidth-2\fboxsep-2\fboxrule\relax}
\begin{tabbing}
+\=++\=++\=++\=++\=++\=++\=++\=\kill
\> FUNCTION $\texttt{minimize}(f,  \bxi, \bar H, R^{\rm init})$ \\
\> $\nu := \Vert \nabla^2f(\bxi)\Vert_1; \quad R :=R^{\rm init}.$ \\
\>  loop \\
\> \> $\Delta \bxi := \texttt{computeDeltaXi}(\nabla^2f(\bxi), \bar H + 10^{-3}\nu I
, \nabla f(\bxi), R).$ \\
\> \> $\bxi^{\rm TEST} := \bxi + \Delta\bxi.$ \\
\> \> if $\bxi^{\rm TEST}$ infeasible \\
\> \> \> $R := R/4.$ \\
\> \> else \\
\> \> \> $\rho:=\frac{f(\bxi) - f(\bxi^{\rm TEST})}{m(\bxi)-m(\bxi^{\rm TEST})}.$ \\
\> \> \> $\rho_g := \frac{\Vert \nabla f(\bxi^{\rm TEST}) - \nabla f(\bxi) -
  \nabla^2 f(\bxi)\Delta \bxi\Vert}{\Vert\nabla f(\bxi)\Vert+
  \Vert\nabla f(\bxi^{\rm TEST})\Vert}.$ \\
\> \> \>  if $\rho < 1/8$ or $\rho_g > 1$ \\
\> \> \>  \>  $R := R/4.$ \\
\> \> \>  else \\
\> \> \> \> if $\rho < 1/4$ \\
\> \> \> \> \> $R := R /2.$ \\
\> \> \> \> else if $\rho \ge 3/4$ and $\lambda > 0$ and $\rho_g\le 0.125$ \\
\> \> \> \> \> $R := 2R$ \\
\> \> \> \> end if \\
\> \> \> \> $\bxi := \bxi^{\rm TEST}.$ \\
\> \> \> end if \\
\> \> end if \\
\> \> if $\lambda\le {\sl tol}_1$ and $\Vert \Delta\bxi\Vert \le {\sl tol}_2$ \\
\> \> \> return $\bxi$, $R$ \\
\> \> end if \\
\> end loop\\
\end{tabbing}
\end{minipage}}
\end{center}
\caption{Minimization algorithm (trust-region method)}
\label{fig:primalonly}
\end{figure}

The routine to compute a single step of the trust-region method
by finding the correct Lagrange multiplier $\lambda$
appears in Fig.~\ref{fig:primalonly_sub} and is standard
(see \cite{NocedalWright}).  The test for positive definiteness as
well as the computation of the direction is carried out with
sparse Cholesky factorization.  The function
$q(\lambda)$ appearing in this procedure is
$$q(\lambda)=\frac{1}{R} - \frac{1}{\Vert N^{1/2}(H+\lambda N)^{-1}\g\Vert};$$
the correct multiplier $\lambda$ should satisfy either $\lambda=0$
or $q(\lambda)=0$, i.e.,
$\Vert N^{1/2}(H+\lambda N)^{-1}\g\Vert = R$.

\begin{figure}[htp]
\begin{center}
\fbox{\begin{minipage}{\dimexpr\textwidth-2\fboxsep-2\fboxrule\relax}
\begin{tabbing}
+\=++\=++\=++\=++\=++\=++\=++\=\kill
\> FUNCTION $\texttt{computeDeltaXi}(H,N,\g, R)$ \\
\> $\lambda_{\rm LOW}:= 0; \quad \lambda_{\rm HIGH}:=\infty; \quad \lambda:=0$ \\
\> loop \\
\> \> if $\lambda_{\rm HIGH}-\lambda_{\rm LOW}<{\sl tol}_3$ \\
\> \> \> Switch to hard-case method (see \cite{NocedalWright}). \\
\> \> end if\\
\> \> $G := H + \lambda N.$ \\
\> \> if $G$ is positive definite \\
\> \> \> $\Delta \bxi := -G^{-1}\g$ \\
\> \> \> $\delta := (\Delta\bxi^TN\Delta\bxi)^{1/2} - R;$ \\
\> \> \> if $|\delta|/R < {\sl tol}_4$ or ($\delta\le 0$ and $\lambda\le{\sl tol}_1$) \\
\> \> \> \> return $\Delta\bxi$ \\
\> \> \> end if \\
\> \> \> if $\delta > 0$ \\
\> \> \> \> $\lambda_{\rm LOW} := \lambda.$ \\
\> \> \> else \\
\> \> \> \> $\lambda_{\rm HIGH} := \lambda.$ \\
\> \> \> end if \\
\> \> \> {\sl PosDef} := {\sl true}; $\lambda := \lambda - q'(\lambda)/q(\lambda).$ \\
\> \> else \\
\> \> \> $\lambda_{\rm LOW} := \lambda.$ \\
\> \> \> {\sl PosDef} := {\sl false}. \\
\> \> end if \\
\> \> if {\sl PosDef} == {\sl false} or $\lambda <\lambda_{\rm LOW}$ or $\lambda>\lambda_{\rm HIGH}$ \\
\> \> \> if $\lambda_{\rm HIGH}==\infty$ and $\lambda == 0$ \\
\> \> \> \> $\lambda := 1$; \\
\> \> \> else if $\lambda_{\rm HIGH}==\infty$ and $\lambda > 0$ \\
\> \> \> \> $\lambda := 2\lambda$; \\
\> \> \> else \\
\> \> \> \> $\lambda := (\lambda_{\rm LOW} +\lambda_{\rm HIGH})/2 $ \\
\> \> \> end if \\
\> \> end if \\
\> end loop \\
\end{tabbing}
\end{minipage}}
\end{center}
\caption{Subroutine of trust-region method to find one step $\Delta\bxi$}
\label{fig:primalonly_sub}
\end{figure}

\section{Computation of energy balance}
\label{sec:energy_balance}

In this section, we describe the terms that enter into the
energy balance used in two of the computational experiments
of Section~\ref{sec:compexp} in order to validate the method.
The energy balances are computed at half-steps between the main
time steps since this is where the displacements are computed
by the implicit midpoint rule.  We note that energy balance
computations for initially rigid cohesive fracture
have been used in the previous literature, e.g., by Molinari et al.\ \cite{Molinari},
to derive results on convergence behavior.

Kinetic energy is evaluated using quadrature of
$$\frac{1}{2}\int_\Omega \rho \dot{u}^2\,dV.$$
The value for $\dot u$ for use in this integral
is taken to be the midpoint of the velocities evaluated at two consecutive
time-steps.  Strain energy is evaluated using quadrature
on the first term of \eqref{eq:twoterm} applied to the displacements
at the midpoint of a time step.

\begin{figure}
  \begin{center}
    \epsfig{file=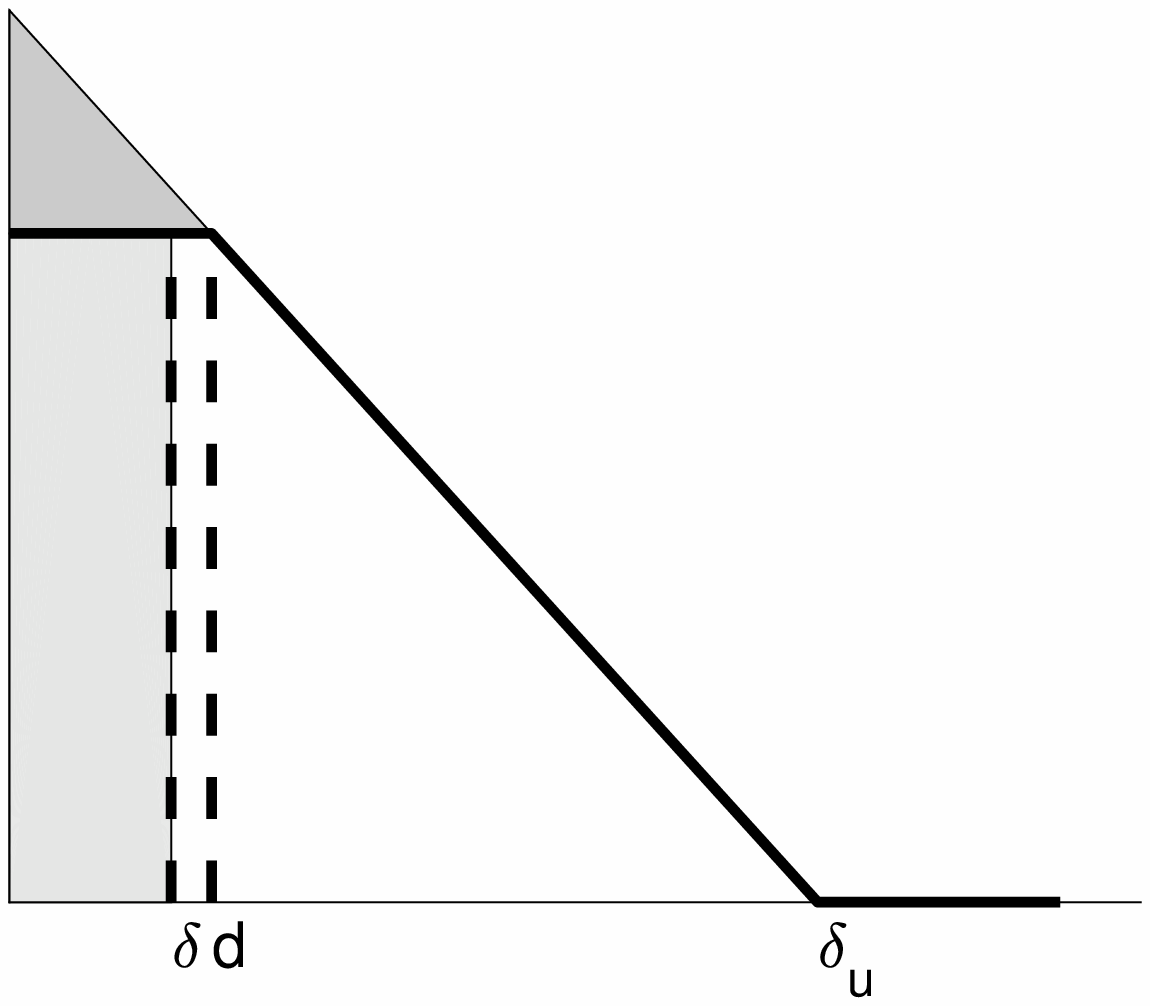,height=4cm} 
    \epsfig{file=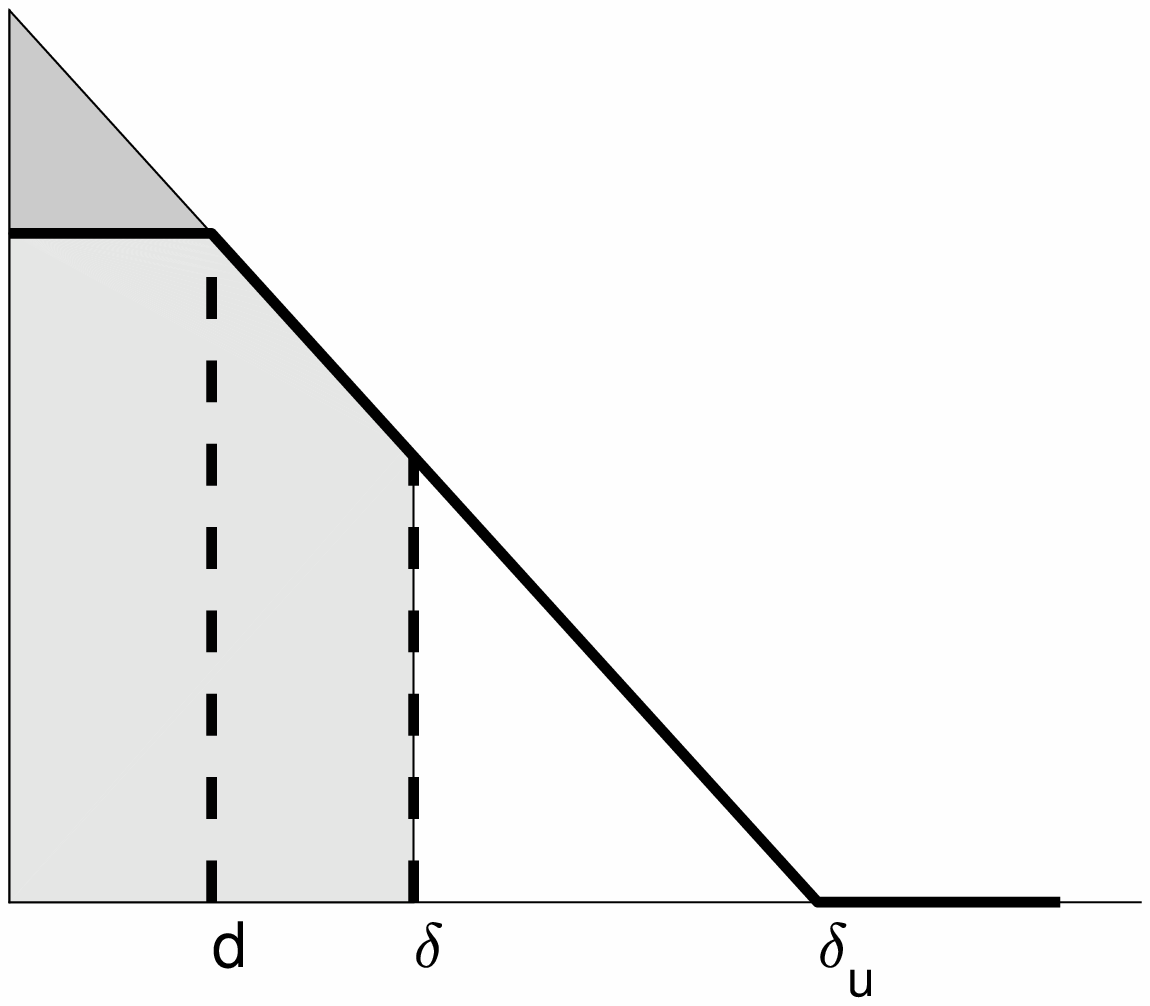,height=4cm} 
  \end{center}
  \caption{Plot of the computation of cohesive energy.  The
    figure in the left considers the case $\delta<d$, where
    $\delta$ is the current opening displacement and $d$ is
    the damage value at an interface Gauss point.  The
    figure on the right considers $\delta>d$.  Note that
    in a continuous-time model, one would never have
    $\delta>d$.  In our computations, however, the update to $d$ lags
    one time-step behind the computation of $\delta$.  The lightly shaded area in both
    figures indicates the energy that can be
    recovered if the interface unloads.  In the plots of the
    next section, this is denoted as ``fracture energy.''  The darkly shaded
    area indicates dissipated energy due
    to irreversibility.
    The total cohesive energy associated with the Gauss
    point is therefore the sum of the two areas.}
  \label{fig:frac_en}
\end{figure}

For a conservative model, cohesive energy would be evaluated using
quadrature on the second term of \eqref{eq:twoterm}.  However, recall that we
have introduced the damage variable $\d$, so the calculation
is more complicated and is described in the caption of
Fig.~\ref{fig:frac_en}.

The stored energies just described must be balanced against the work
done on the models.  The first source of work is from the traction
and body forces (the last two terms of \eqref{eq:twoterm}).
Since our examples do not involve either traction or
body forces,  we omit a detailed discussion of these terms,
but their computation is relatively straightforward.

The work done by
moving displacement (i.e., velocity) boundary conditions 
is computed as follows.  The objective
function $f(\u;\s_0)$ is written down as in \eqref{eq:primbar2},
except that all functions are written in terms of $\u$
(the vector of all DOFs) instead of $\x$ (the vector
of unconstrained DOFs).  Next, the gradient with respect
to $\u$ is computed at the minimizer for the
final value of $\mu$.
The gradient entries corresponding to the
constrained DOFs will in general not vanish because
the objective function is not minimized with respect
to them.  In fact, the gradient entries are exactly
the reaction forces for those DOFs.
Therefore, the work due to displacement constraints during
a single time-step, that is, from one half-time-step to the next,
is evaluated as follows.  One computes the inner
product of the time-average of these forces (average between
the current half-time-step and the previous) and the distance
traveled by each such DOF between the current and
previous half-time-step.

Finally, work done by contact boundary conditions is
also obtained as an inner product.
One multiplies the contribution to the force balance
from the derivative of the contact barrier term
that appears in $\nabla_\u f(\u;\s_0)$ by the distance
traveled.  Time-averages are used as in the last paragraph.
Note that although the coefficient $\mu$ in front of the
barrier term may vanish, the corresponding term of  the
gradient does not vanish as $\mu\rightarrow 0$ but instead
tends to a constant value; this is a well-known aspect
of interior-point theory.

\section{Computational experiments}
\label{sec:compexp}
In this section we describe three computational experiments.
Experiment 1 involves
impact of a metal striker on a compact compression specimen
(CCS).  This application demonstrates the ease in which convex
constraints can be included in the computation.  This problem
involves three contact surfaces detailed in the next paragraph.

\begin{figure}
  \vspace{0.1in}
  \begin{center}
    \epsfig{file=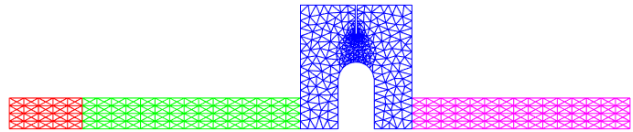, height=8cm} \\
    \epsfig{file=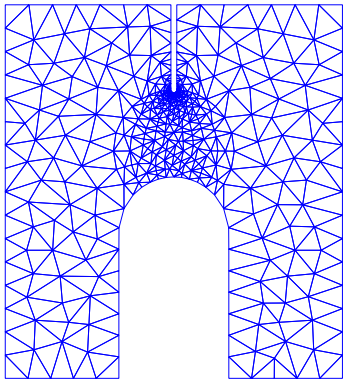, height=5cm}
  \end{center}
   \vspace{0.1in}
   \caption{Initial geometry and mesh of the striker and CCS specimen
     (Experiment 1).
     The entire mesh is shown on top with each
     part (striker, left bar, CCS, right bar) in a different
     color; the initial gap between
     the striker and left bar is not discernable at this
     scale.  A close-up of the CCS initial mesh appears on the bottom.}
  \label{fig:ccsorig}
\end{figure}

The CCS is made of
PMMA of $\mbox{height}\times\mbox{width} = 51\mbox{mm}\times 46\mbox{mm}$.
Its initially undeformed configuration is depicted in 
Fig.~\ref{fig:ccsorig}.  It is initially in contact with two steel bars
(left and right bar), both of  which are
stationary.  Before contact, the striker
is moving in the positive $x$ direction at 25 m/s, and 
a gap exists between the striker and left bar.
This computation simulates an experiment
by Rittel and Maigre \cite{MRC:23:475}.
The striker and the left bar,
the left bar and the specimen, and finally the specimen and
the right bar are all contact surfaces.
We model contact in each of the three surfaces
via inequalities between $x$-coordinates of matching
nodes of the two sides of the surface.  Refer to Fig.~\ref{fig:ccsorig}.
In our formalism, these inequalities are presented as $E\x\ge\a$
in \eqref{eq:mainprob}.  Note that prevention of
interpenetration between neighboring bulk elements is 
modeled as a different set of  linear inequalities in \eqref{eq:mainprob}, namely,
the inequalities $(s_1)_{e,\iota}\ge 0$.

No boundary conditions are applied, i.e., all boundaries are
unconstrained and traction-free.
The initial mesh of the CCS contains $1765$ nodes 
and $820$ element (quadratic triangles).  After duplication
of nodes  to create interface elements, the number of
nodes is $6\cdot 820=4920$. 
The mesh near the fracture zone
is an isoperimetric ``pinwheel'' mesh \cite{pinwheel1} transformed by
a nonlinear coordinate transformation, as depicted in
Figure~\ref{fig:ccsorig}.  Isoperimetric meshes have the property
that in the limit of mesh refinement, all possible crack orientations
are represented in the mesh, so they are well suited for computations
in which determining the crack path is part of the problem.
(Other techniques have been proposed
for representing many possible orientations
in a mesh, e.g., adaptive splitting of polygonal elements
by Leon et al.\ \cite{Paulino}.)
The three metal bars (striker, left bar, and right bar)
are modeled with quadratic triangles.
No interface
elements are introduced in the three steel bars.  
The final model
contains 6137 nodes and 1380 quadratic triangles.

The  hyperelastic model
used in the CCS is nonlinear
plane stress,
whose material parameters are $c_1$ and $\beta$ (see \eqref{eq:knowles}).
These are obtained from the reported
$E$ and $\nu$ according to the formulas $c_1=E/(4(1+\nu))$,
$\beta=\nu/(1-2\nu)$ \cite{KnowSter1983}.
The striker is isotropic linearly elastic.

\begin{figure}
  \vspace{0.1in}
  \begin{center}
    \epsfig{file=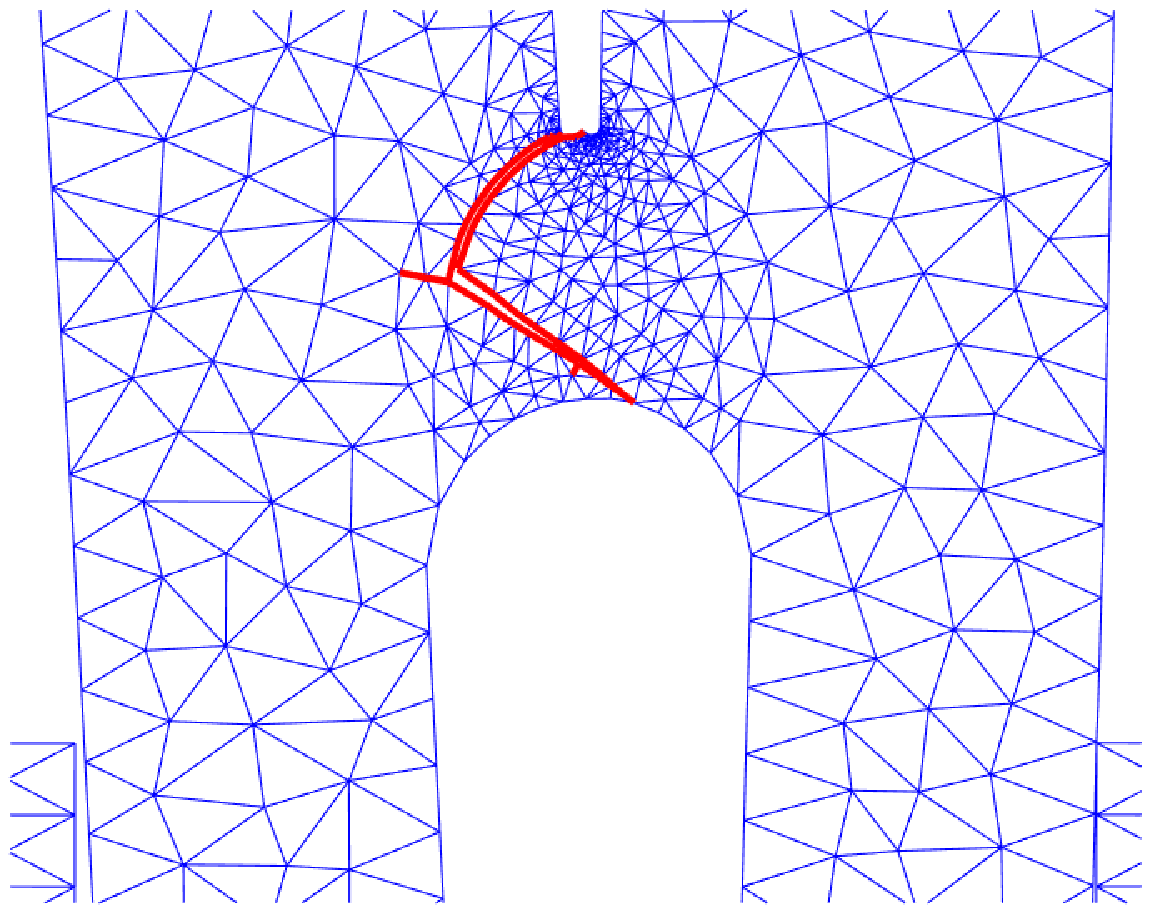, height=5cm}
  \end{center}
  \vspace{0.1in}
  \caption{Damaged interfaces after 40 time steps (120 microseconds
    of simulated time) in Experiment 1.}
  \label{fig:ccsredlines}
\end{figure}

The time step is $3$ microseconds.  The simulated crack path after 40
time steps (120 simulated microseconds)---refer to
Fig.~\ref{fig:ccsredlines}---is similar
to experimental results in \cite{MRC:23:475}.

The test was run on a 2.6GHz Intel Xeon E5-2690 running Linux.
The algorithm was coded in the Julia \cite{Julia2017} programming
language, version 1.1.0.  The core computing kernel is sparse
Cholesky factorization, which uses Julia's implementation of
SuiteSparse \cite{DavisBook}.
The code is not parallelized yet.  The computation time for 40 time
steps was 2 hours.

\begin{table}[htp]
\begin{center}
\begin{tabular}{llrrrr}
\hline
Property & Symbol & PMMA & Steel  & Concrete & Mortar  \\
         &        & (Exper.~1) &(Exper.~1)&(Exper.~2) & (Exper.~3) \\
\hline
Young modulus & $E$ & 5.76 GPa  & 200 GPa & 38 GPa & 5.98 GPa\\
Poisson ratio & $\nu$ & 0.42 & 0.3  & 0.18 & 0.22\\
Density & $\rho$ & 1180 kg/m$\mbox{}^3$ & 8050 kg/m$\mbox{}^3$ &--- & ---\\
Critical traction & $\sigma_c$ &105 MPa & --- & 3 MPa & 3 MPa\\
Mixity  & $\beta^{\rm MIX}$ & 2.0 &  --- & 1.5 & 1.0 \\
Critical energy \\
release rate & $G_c$ & 352 Pa$\cdot$m &  --- & 69 Pa$\cdot$m & 2280 Pa$\cdot$m\\
\hline
\end{tabular}
\end{center}
\caption{Material properties used in computational experiments.}
\label{tab:matprop}
\end{table}

We also tracked energy balance. The technique used
to measure energy balance was described in Section~\ref{sec:energy_balance}.
Because the optimization problem is solved at the
midpoint of timesteps in the implicit midpoint rule used herein, we
evaluate the energy balance at time-step midpoints.
However, not all the variables are evaluated at time-step midpoints, so interpolations
must be used.  Therefore, we would not expect exact energy balance because the
different terms involve different approximation assumptions.

\begin{figure}
  \begin{center}
    {\epsfig{file=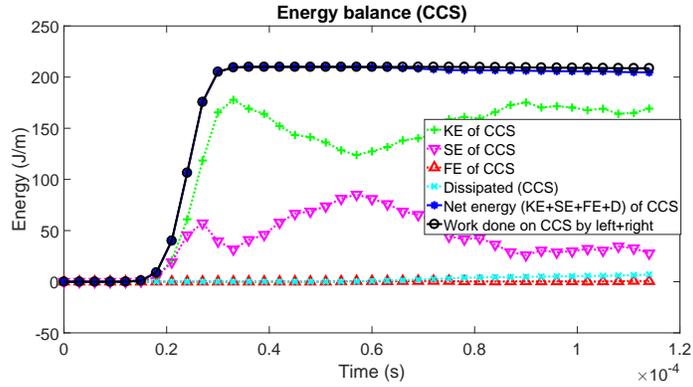,height=5cm}} \\
    (a) Energy balance of CCS \\
    {\epsfig{file=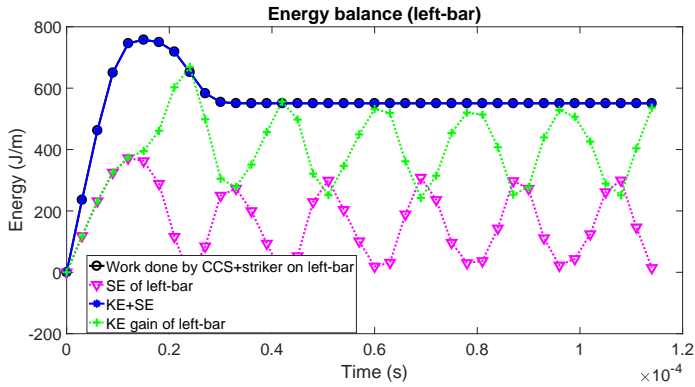,height=5cm}} \\
    (b) Energy balance of left bar
  \end{center}
  \caption{Energy balances for the CCS (Experiment 1); ``KE'' stands for
    kinetic energy, ``SE'' stands for strain energy, and
    ``FE'' stands for fracture energy.  Refer
    to Section~\ref{sec:energy_balance} for an explanation of
    how the energy and work contributions were computed.
    In the top figure, if energy were exactly
    conserved, the blue and black curves would coincide.
    In the bottom figure, the curves do coincide.   The
    energy balance for the left bar indicates that it
    vibrates after the impact.}
  \label{fig:ccs_energy}
\end{figure}

\begin{figure}
  \begin{center}
    \epsfig{file=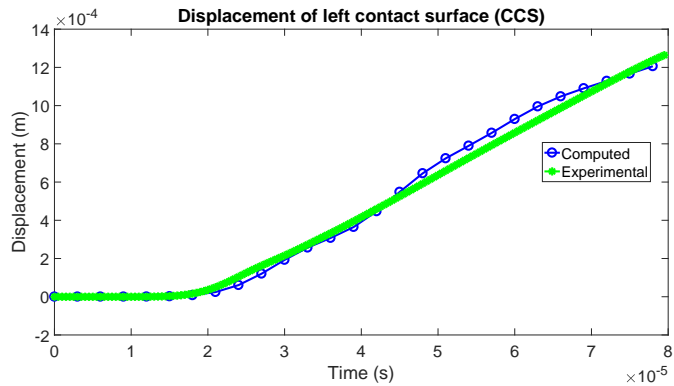, height=5cm}
  \end{center}
  \caption{Comparison of experimentally measured displacement of the
    CCS left contact surface with the displacement computed by the
    method (Experiment 1).  The latter is averaged over the contact surface.}
  \label{fig:contactsurfdispl}
\end{figure}

Fig.~\ref{fig:contactsurfdispl}
compares the horizontal displacement of the left contact
surface of the CCS from the computation versus an experimental
measurement \cite{rittelpriv}.
The comparison is imprecise because only a single
number was recorded per time step by the experiment. On the other hand, the
computation indicates that the left contact surface of the CCS bends
inward, and therefore its horizontal displacement depends on
the vertical coordinate of the measurement.
To account for this, we compare the
average computed horizontal displacement over the entire
contact surface to the
experimental measurement.  Note that the time-coordinate of the
experiment was shifted by hand to match the start time of the
computation since the experimental data did not identify
the time value when the striker collides with the left bar.

\begin{figure}
  \vspace{0.1in}
  \begin{center}
    \epsfig{file=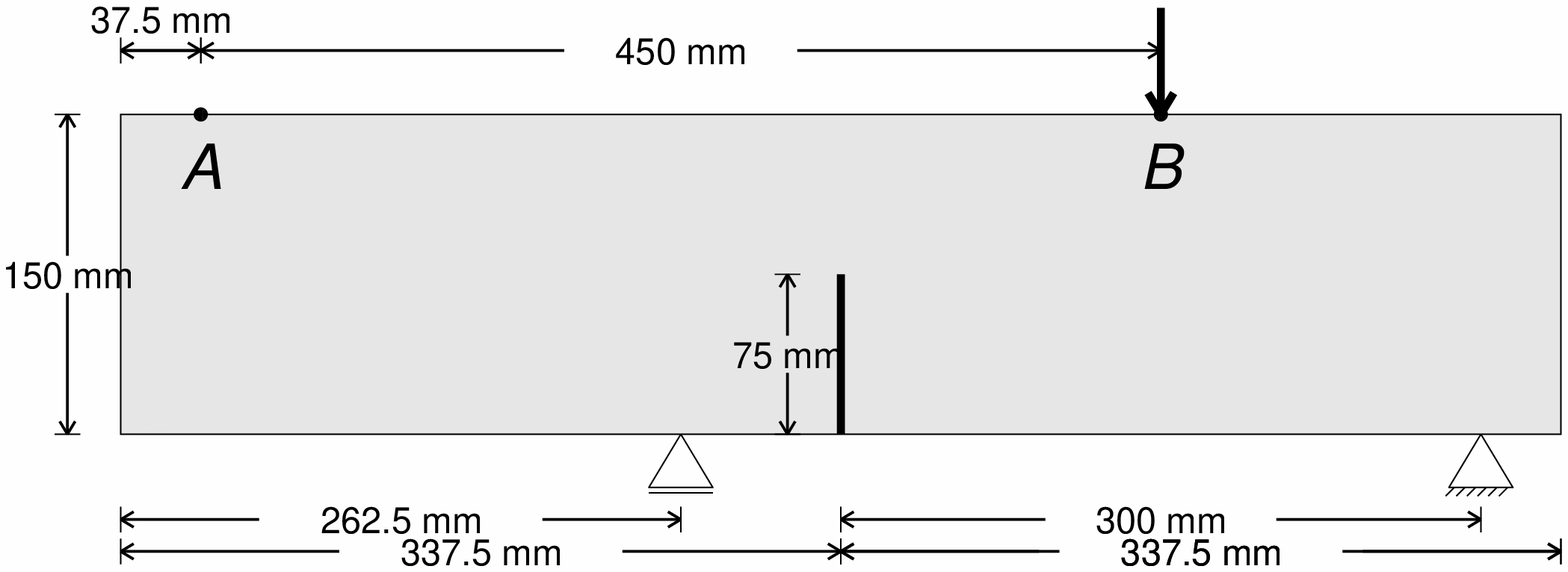, height=3.8cm}
  \end{center}
  \vspace{0.1in}
  \caption{Experiment 2: quasistatic test of concrete beam.  The beam is supported
    at two asymmetric points of its base and has an initial vertical slit.
    The point B is loaded quasistatically.  The point A is referred to
    in Fig.~\ref{fig:galv_load_defl} below.}
  \label{fig:galv1}
\end{figure}

Experiment 2 exhibits the performance of the method
in the quasistatic (slow loading) regime.  In this problem, taken
from Galvez et al.\ \cite{Galvez}, a beam of
concrete 15 cm $\times$ 67.5 cm, plane stress, with an initial
centered slit is 
subjected to a moving displacement boundary condition
concentrated at an off-center point.  The problem set-up is
described in more detail in Fig.~\ref{fig:galv1}.

\begin{figure}
  \begin{center}
    \epsfig{file=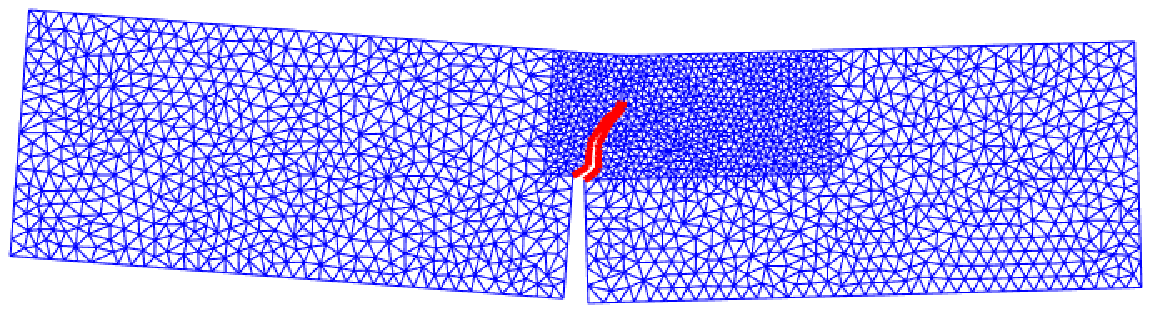, height=7.0cm}
  \end{center}
  \caption{Final configuration of concrete beam after 25 steps.
    The displacements are exaggerated by a factor of 50
    for better visualization.}
  \label{fig:galv_final}
\end{figure}

A depiction of the configuration after 26 load steps is shown in
Fig.~\ref{fig:galv_final}.  The number of elements is
$3168$ and the number of nodes (after duplication)
is $19,008$.  This computation required 17 hours.
(The amount would be greatly reduced if we had inserted interface
elements only in the zone where crack propagation is known to occur i.e.,
above and to the right of the initial slit, whereas in fact our mesh
has cohesive interfaces at every interelement boundary.)
The crack path roughly matches the experiment in \cite{Galvez},
although we were not trying to accurately reproduce the path in this
experiment because we did not use a pinwheel or other special mesh.

\begin{figure}
  \begin{center}
    \epsfig{file=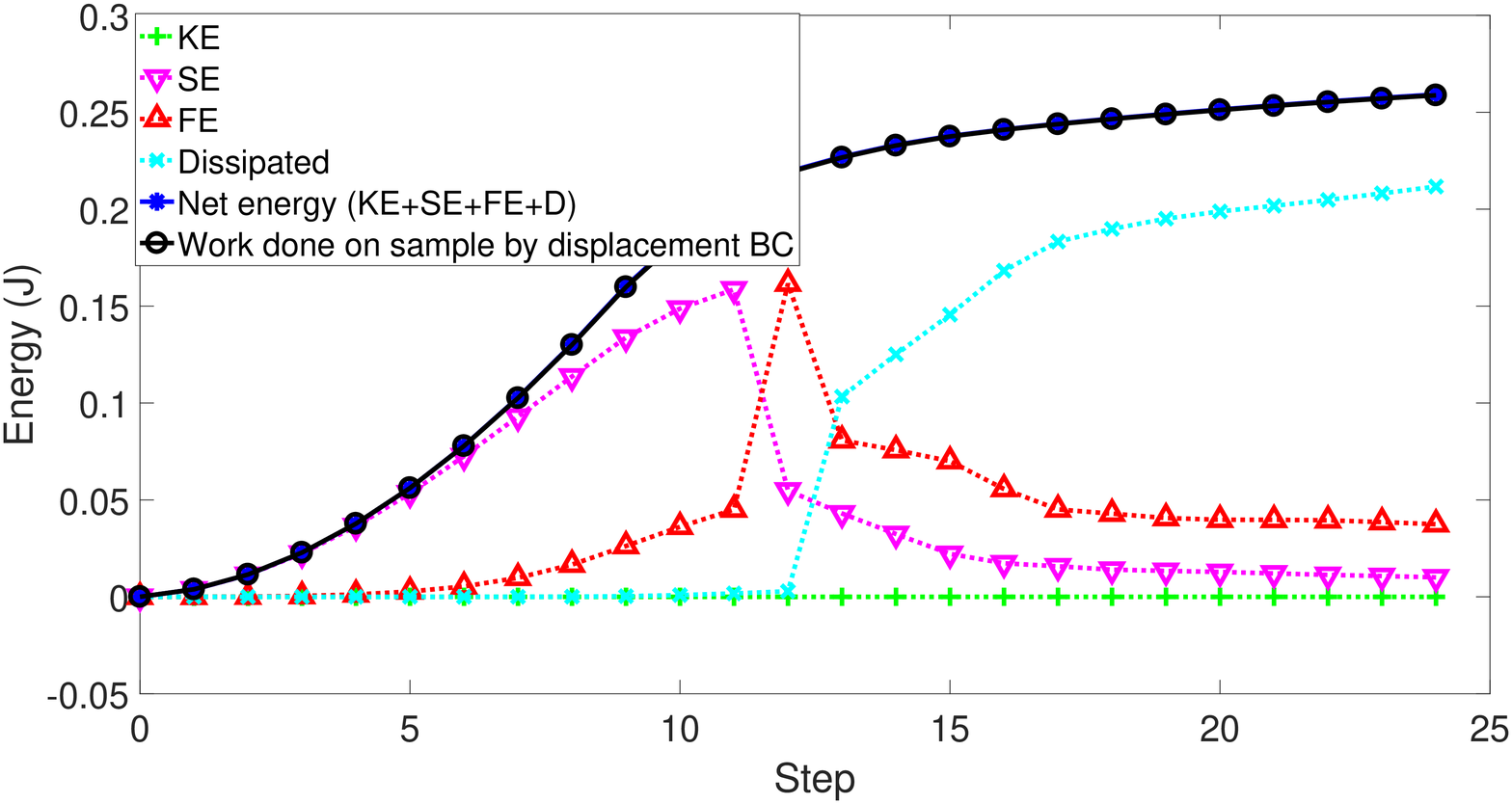, height=7.0cm}
  \end{center}
  \caption{Energy balance for Experiment 2.  Refer to the caption
    of Fig.~\ref{fig:ccs_energy} for further information. 
    The last two curves are indistinguishable (coincident).  Thickness has
    been normalized to 5cm to correspond to the experiment.}
  \label{fig:galv_energybalance}
\end{figure}

An energy balance was also computed for the Galvez experiment; the
results are reported in Fig.~\ref{fig:galv_energybalance}.
It is interesting to compare the energy balance in Experiment 1 versus
Experiment 2.   In Experiment 1, most of
the work on the CCS becomes kinetic energy of the CCS.  Part
becomes strain energy, but the strain energy is partly released as it is
transformed first to cohesive energy and then dissipated.  In Experiment 2,
there is no kinetic energy (the problem is
quasistatic), so the energy of the load first builds up the strain
energy, which is then released as cohesive energy, and then soon after
the fracture energy is dissipated.

\begin{figure}
  \begin{center}
    \epsfig{file=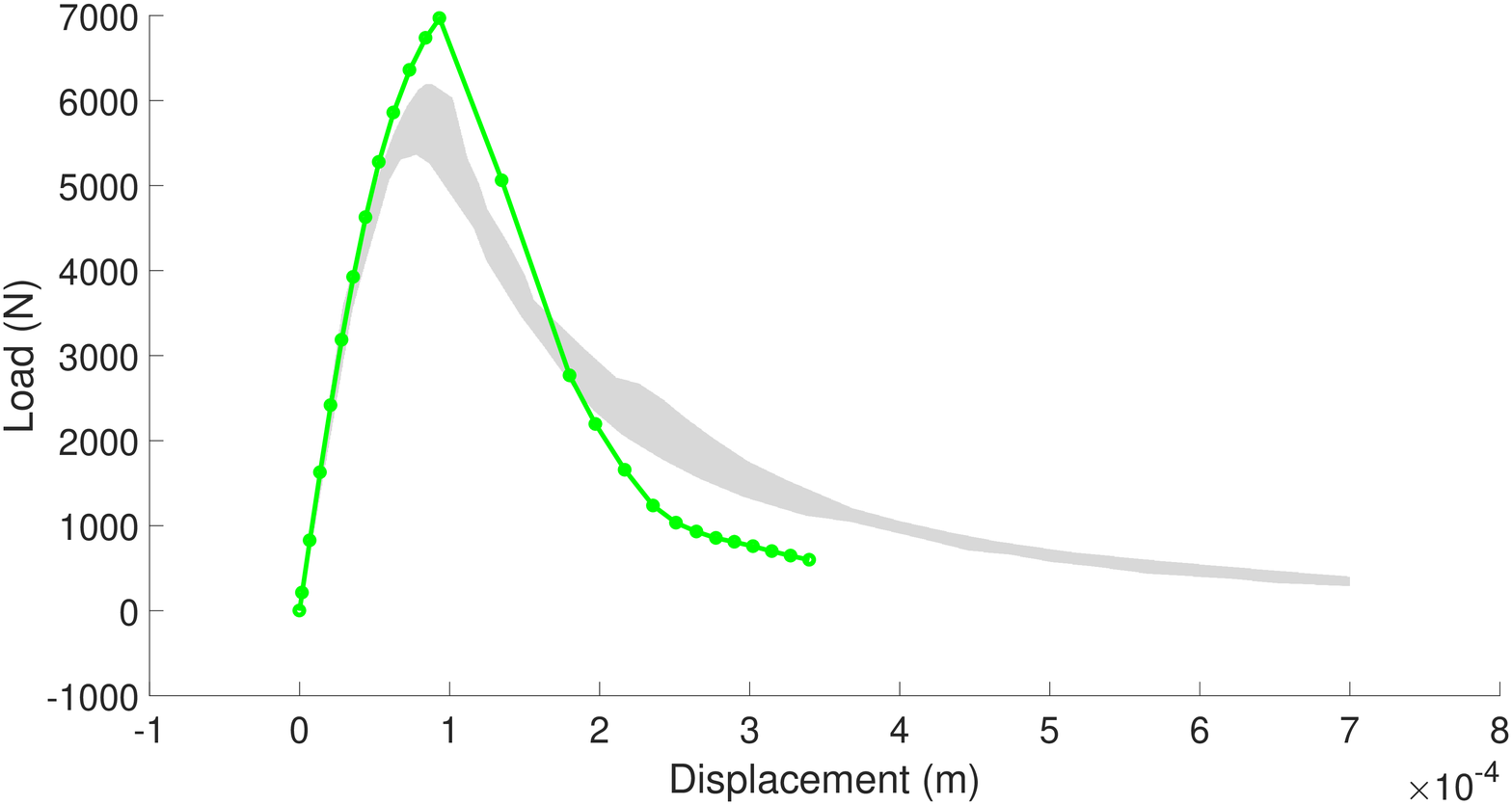, height=5.0cm}
  \end{center}
  \caption{A plot of the load (reaction force) at point B
    versus the vertical deflection of point A (experiment 2).
    The green curve represents the data from our computation.  The experimental
    envelope from Galvez et al.\ \cite{Galvez} is represented
    by gray curves.  Thickness has
    been normalized to 5cm to correspond to the experiment.}
  \label{fig:galv_load_defl}
\end{figure}

In addition, a load-displacement curve was plotted since this data
is available experimentally.  The result of this computation
appears in Fig.~\ref{fig:galv_load_defl}.

\begin{figure}
  \begin{center}
    \epsfig{file=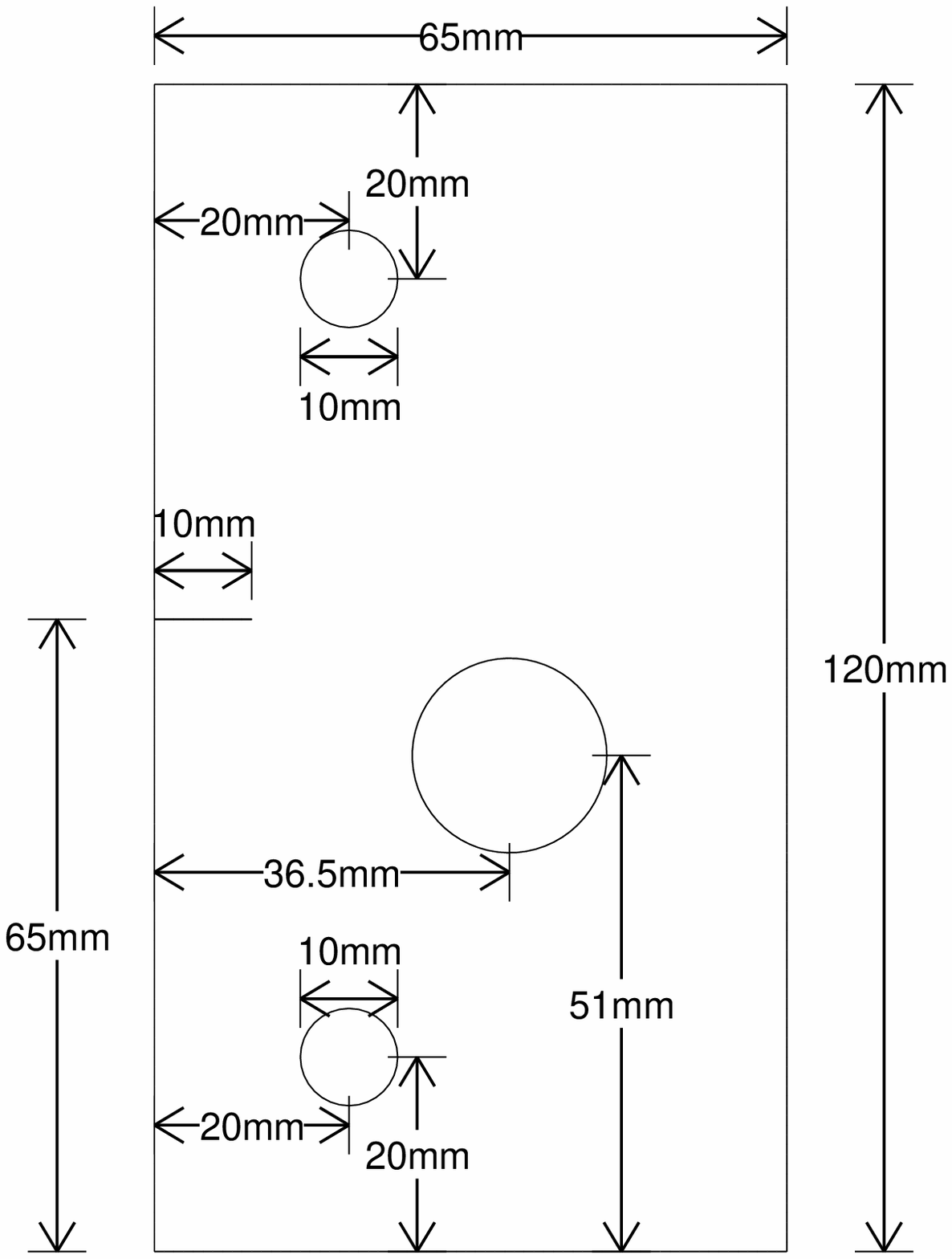, height=7.0cm}
  \end{center}
  \caption{Experiment 3 involves
    a 65mm$\times$120mm concrete plate with three holes and a 10mm-long
    horizontal notch on the left edge.  The plate is loaded quasistatically
    by pins that pull the top and bottom holes vertically apart.}
  \label{fig:threeholediagram}
\end{figure}

Experiment 3 is
also a quasistatic experiment involving a concrete mortar plate with three
holes and a notch illustrated in Fig.~\ref{fig:threeholediagram}.
Pins that fit into the top and bottom holes 
pull the holes vertically apart at a rate
$0.1$ mm/s.  This is modeled as velocity boundary conditions constraining
both $x$- and $y$-coordinates of all points on the boundaries of these two holes.

Ambati et al.\ \cite{Ambati2015}
present a computational result using their phase-field method
to determine the load-displacement relationship (displacement of the top
pin versus force on the top pin).  Their result shows two peaks in
the load, the first at slightly more than 0.5 mm  and the second at
slightly less than 2 mm.

\begin{figure}
  \begin{center}
    \epsfig{file=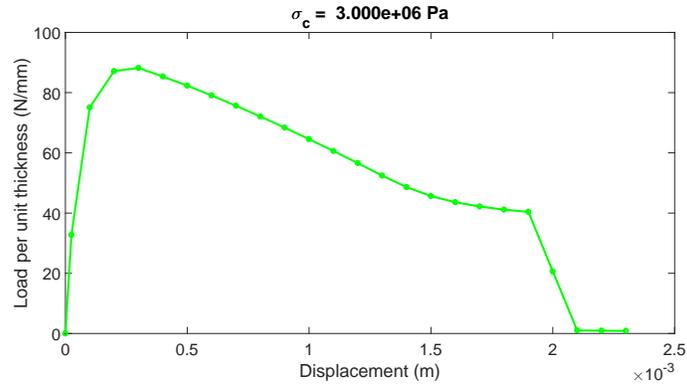, height=5.0cm} \\
    (a) \\
    \epsfig{file=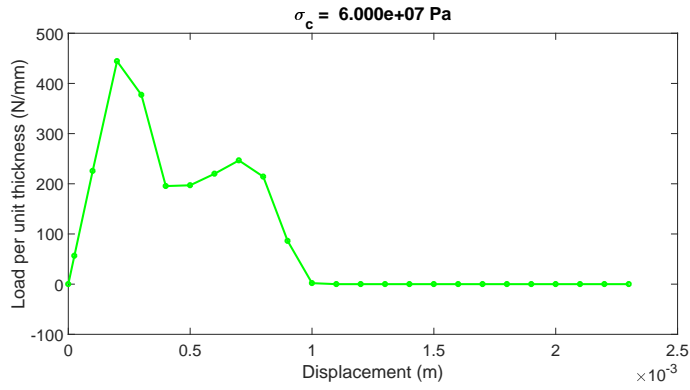, height=5.0cm} \\
    (b) \\
    \epsfig{file=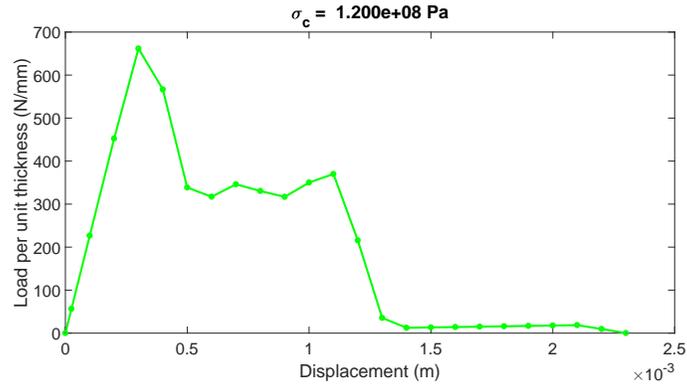, height=5.0cm} \\
    (c) \\
  \end{center}
  \caption{Load-displacement curve for Experiment 3 with three
    different values of the material parameter $\sigma_c$.
    Higher values of $\sigma_c$ are able to reproduce
    the double peaks reported in Ambati et al. \cite{Ambati2015}.}
  \label{fig:threehole_load_disp}
\end{figure}

Our load-displacement plot is presented in Fig.~\ref{fig:threehole_load_disp}.
This computation required six hours per parameter choice. The load was applied at an incremental
rate of 0.1 mm per step.
We were
unable to reproduce the second peak with the material parameters shown in
Table~\ref{tab:matprop}.  We hypothesize the following
phenomenological explanation for the second peak.
After the crack propagates to the large central hole, the body
behaves elastically (i.e., the remaining reverse C-shaped piece
flexes elastically) until a second crack starts on the right of the big hole.
Our computational experiment, however, shows that the second crack on the
right of the big hole is mostly formed by time the first crack reaches
the big hole, so no elastic behavior is observed after the first crack
reaches the big hole.

The phenomenological explanation in the previous paragraph suggests that
our method could obtain the second peak if we increase $\sigma_c$,
thus delaying the nucleation of the  second crack.  We tested this
hypothesis by raising $\sigma_c$  by a factor of 20
(to $6.0\cdot 10^7$ Pa) and by a factor of  40 (to $1.2\cdot 10^8$ Pa).  Note
that these numbers are significantly higher than the usual reported
critical stress for concrete mortar.  With these modified values
of $\sigma_c$ we indeed observed a second peak as in Fig.~\ref{fig:threehole_load_disp},
but even at
$\sigma_c=1.2\cdot 10^8$ the displacement at the position
of the second peak falls well short of 2mm.  

The phase-field method used in the Ambati computation does not have
$\sigma_c$ as a material parameter.
The discussion in the previous paragraph raises the question of what value of
$\sigma_c$ corresponds to the Ambati
computation. 
Some authors (see, e.g., eq. (27) of Borden et al.~\cite{borden2012}) have proposed
that the length-regularization parameter in phase-field models is linked to $\sigma_c$.
Experiment 3 may provide an
example for future investigation of the connection between length-regularization
in such phase-field methods and $\sigma_c$ in cohesive models. A phase-field
regularized cohesive model that involves $\sigma_c$ as a material parameter and
frees previous implementations of phase-field method from the above limitation
has been recently proposed in Geelen et al.~\cite{Dolbow2}.

\begin{figure}
  $$
  \begin{array}{ccc}
    \epsfig{file=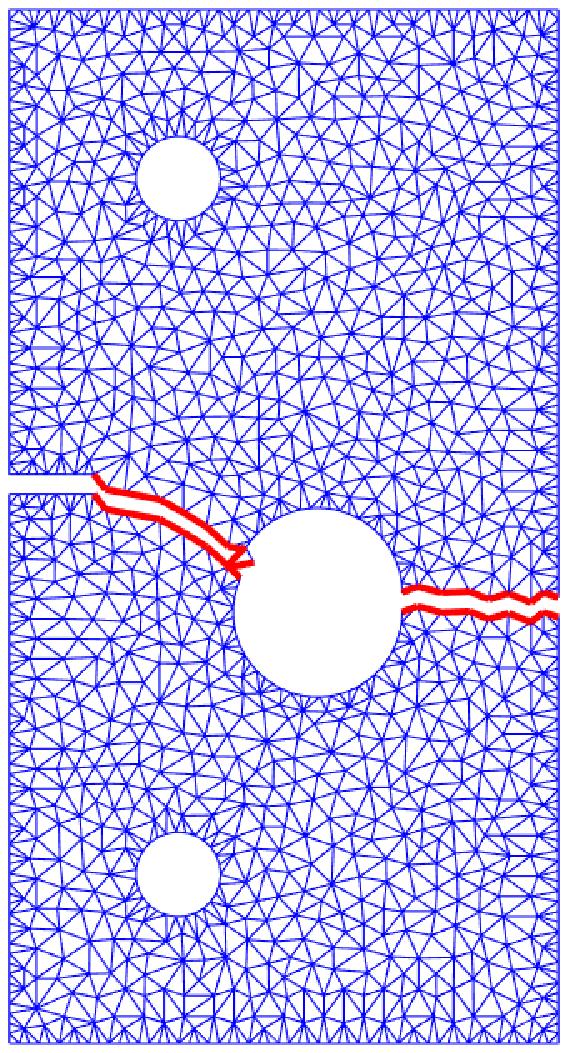, height=3.5cm} & 
    \epsfig{file=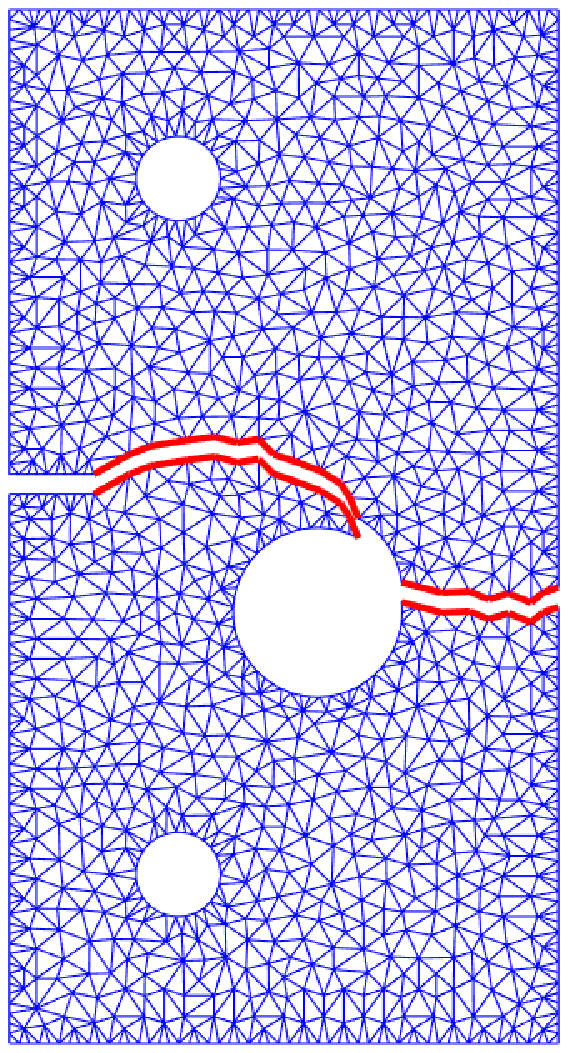, height=3.5cm} & 
    \epsfig{file=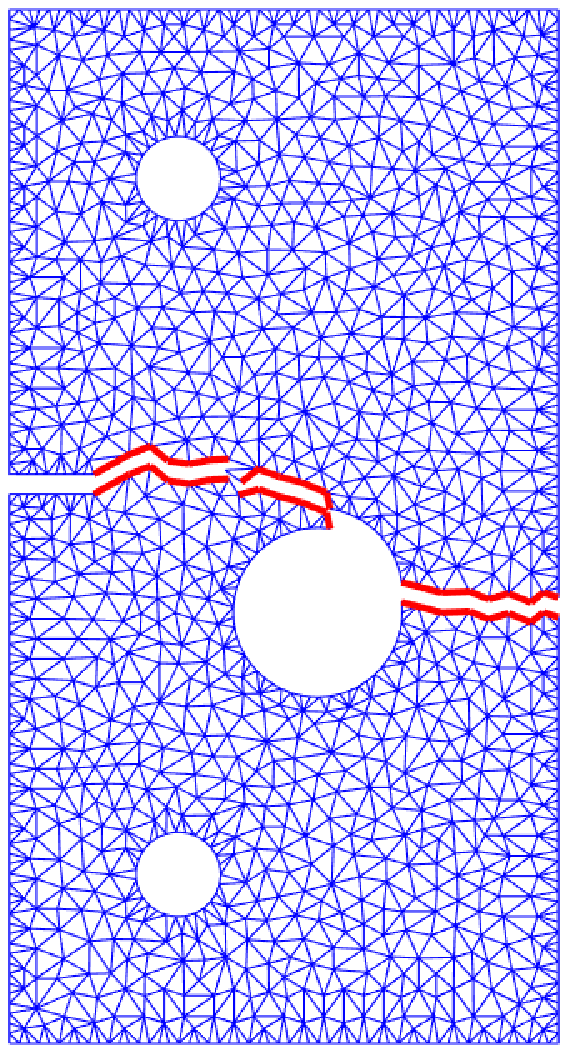, height=3.5cm}
    \\
    \mbox{(a)} & \mbox{(b)} & \mbox{(c)}
  \end{array}
  $$
  \caption{Crack-path in Experiment 3 for three different values of
    $\sigma_c$.}
  \label{fig:threehole_crackpath}
\end{figure}

The crack-paths for all three values of $\sigma_c$ are shown Fig.~\ref{fig:threehole_crackpath} as well as the laboratory experimental path.
The first part of the crack path is sensitive to $\sigma_c$ but the second
part is not.

A fourth computational experiment not reported here
was a branching dynamic crack model
used frequently in the literature (see, e.g., Borden et al.\ \cite{borden2012})
in which a rectangular specimen with a horizontal
notch is subject to vertical traction loading.  This traction load causes a horizontal
crack to propagate from the notch.
After a small distance
of horizontal propagation, the crack
branches, yielding a crack with a tilted ``Y'' shape.
We found that our method produced
a branch with the correct
shape, but the straight part of the crack was too long compared to other computations
in the literature.  We subsequently
observed that the branch started at the correct position
when the method was re-run with a very small time step,
requiring a simulation that ran for nearly
a week.
This is a limitation not of our method specifically but rather of all methods based
on implicit time-step rules.  In problems in which the time step must be small
enough to capture highly dynamic and unstable behavior, implicit
time-stepping is generally not recommended.  Refer to the recent
paper \cite{HirmandPapouliaBlockCD} by the second and third authors
for an explicit time-stepping method that uses energy principles similar
to those herein.

\section{Summary and conclusions}
\label{sec:conc}
An interior-point method for initially rigid cohesive fracture
is proposed. A key technical step to make this method possible
is the replacement of an equation relating the effective opening
displacement to the coordinate entries of the opening displacement
by an inequality constraint.

The specific optimization models introduced herein are as follows.
\begin{itemize}
\item
  Model \eqref{eq:optmodeldisc} is the fundamental optimization
  problem that introduces all the
  energy terms into the objective and the constraints relating
  the bulk-node displacements $\x$ to opening displacements $\s_0$.

\item
  Model \eqref{eq:mainprob} is identical to 
  \eqref{eq:optmodeldisc} except for replacing an equality
  constraint with an inequality constraint that yields
  a mathematically equivalent optimization problem (i.e., same optimizers).
  This replacement makes the model amenable to an interior-point method.
\item
  Model \eqref{eq:primbar} replaces the conic inequality constraints
  in \eqref{eq:mainprob} with self-concordant barrier functions in
  the objective
  as is the usual practice in the development of an interior-point 
  method.  As the barrier parameter $\mu>0$ tends to 0, a solution
  to \eqref{eq:mainprob} is recovered from the solution to \eqref{eq:primbar}.
\item  
  Model \eqref{eq:primbar2} is mathematically equivalent to
  \eqref{eq:primbar}; the modification is that substitution has
  been used to eliminate the equality constraints.  The main
  work of the code is the solution of \eqref{eq:primbar2}.
\item
  Model \eqref{eq:primbar2phaseI0} is the artificial problem
  solved in Phase I to obtain an initial feasible solution to
  \eqref{eq:primbar2}.  This model introduces the artificial
  variable $t$.
\item
  Model \eqref{eq:primbar2phaseI} further develops
  \eqref{eq:primbar2phaseI0} to obtain an initial feasible
  solution in which all elements are non-inverted even in the
  presence of moving displacement boundary conditions.  This
  is accomplished by changing variables from $\x$ (unconstrained
  boundary DOFs) to $\u$ (all boundary DOFs).  The boundary
  conditions are represented by equality constraints on $\u$, but
  \eqref{eq:primbar2phaseI} enforces  these equality constraints
  as inequalities involving the artificial variable $t$ to
  ensure that a feasible starting point exists for a known
  solution $\u$ with no inverted elements.
\end{itemize}

A summary of the computational procedures presented is as follows.
\begin{itemize}
\item
  The top-level procedure appears in Fig.~\ref{fig:toplevel}.  Its loop
  is over time-steps or load steps.  It also maintains and updates
  the damage variables.  In addition, on every third
  step it recomputes the matrices used in regularizing the
  trust-region method.
\item
  The two main loops of the interior-point method
  appear in Fig.~\ref{fig:solvermainprob}.  The first loop
  is Phase I in which \eqref{eq:primbar2phaseI0}  is solved
  in order to obtain a feasible starting point for
  \eqref{eq:primbar2}.  The second loop is Phase II
  in which \eqref{eq:primbar2} is solved for a decreasing
  sequence of $\mu$'s.
\item
  The solver for a specific instance of either
  \eqref{eq:primbar2phaseI0} or \eqref{eq:primbar2}
  (i.e., for one particular value of $M$ or $\mu$)
  is a trust-region method, which is
  detailed in Fig.~\ref{fig:primalonly} and 
  Fig.~\ref{fig:primalonly_sub}.  The trust-region
  method is necessary due to the nonconvexity of the
  energy functional; it replaces the Newton loop that would
  be present in a conventional convex interior
  point method.
\end{itemize}

Computational tests show that
the method is practical for quasistatic and moderately fast
dynamic problems and can easily encompass additional
conic inequality constraints.

\section{Acknowledgements}
The authors thank M.~Ambati for helpful discussion regarding Experiment 3.
The authors thank the referees for their careful reading and valuable
suggestions for revising the paper.

\bibliographystyle{plain}
\bibliography{mybib,mypapers,intpt}
\end{document}